\documentclass[11pt]{amsart}
\usepackage[T1]{fontenc}
\usepackage{lmodern,amsmath,amssymb,amsthm,mathtools,mathrsfs}
\usepackage[letterpaper,margin=0.97in,headheight=14pt]{geometry}
\usepackage{microtype,booktabs,array,longtable,enumitem}
\usepackage{needspace}
\usepackage{tikz-cd}
\usepackage[unicode,colorlinks=true,linkcolor=black,citecolor=black,urlcolor=blue]{hyperref}
\hypersetup{
  pdftitle={Algebraically elliptic resolutions in a family of determinantal quartic threefolds},
  pdfauthor={Alexander Dvorsky, Shulim Kaliman, and Mikhail Zaidenberg}
}
\usepackage{bookmark}
\setlist[enumerate]{topsep=4pt,itemsep=2pt,parsep=0pt,leftmargin=2em}
\newtheorem{theorem}{Theorem}[section]
\newtheorem{proposition}[theorem]{Proposition}
\newtheorem{lemma}[theorem]{Lemma}
\newtheorem{corollary}[theorem]{Corollary}
\theoremstyle{definition}
\newtheorem{definition}[theorem]{Definition}
\newtheorem{notation}[theorem]{Notation}
\newtheorem{remark}[theorem]{Remark}
\numberwithin{equation}{section}
\newcommand{\PP}{\mathbf P}
\renewcommand{\AA}{\mathbf A}
\newcommand{\CC}{\mathbf C}
\newcommand{\QQ}{\mathbf Q}
\newcommand{\ZZ}{\mathbf Z}
\newcommand{\OO}{\mathcal O}
\newcommand{\Bl}{\operatorname{Bl}}
\newcommand{\Sing}{\operatorname{Sing}}
\newcommand{\Dom}{\operatorname{Dom}}
\newcommand{\Indet}{\operatorname{Indet}}
\newcommand{\rk}{\operatorname{rank}}
\newcommand{\Pic}{\operatorname{Pic}}

\newcommand{\sm}{\mathrm{sm}}
\newcommand{\red}{\mathrm{red}}
\newcommand{\cS}{\mathcal S}\newcommand{\cC}{\mathcal C}
\newcommand{\cL}{\mathcal L}

\newcommand{\coeffmatrix}{\begin{pmatrix}-2&1&0&-2&3\\-3&0&1&1&1\\0&-3&2&-1&2\\-1&-2&0&2&1\end{pmatrix}}
\title[Elliptic resolutions of determinantal quartics]{Algebraically elliptic resolutions\\in a family of determinantal quartic threefolds}
\author{Alexander Dvorsky}
\address{Department of Mathematics, University of Miami,
Coral Gables, FL 33124, USA}
\email{dvorsky@gmail.com}
\author{Shulim Kaliman}
\address{Department of Mathematics, University of Miami,
Coral Gables, FL 33124, USA}
\email{kaliman@math.miami.edu}
\author{Mikhail Zaidenberg}
\address{Univ. Grenoble Alpes, CNRS, IF,
38000 Grenoble, France}
\email{mikhail.zaidenberg@gmail.com}
\date{}

\expandafter\providecommand\csname subjclassname@2020\endcsname{2020 Mathematics Subject Classification}

\begin{document}
\begin{abstract}
We construct a family of determinantal irrational quartic threefolds in $\PP^4$ for which the blow-up of the sixteen nodes and two distinguished small projective resolutions are algebraically elliptic in the sense of Gromov.
\end{abstract}
\maketitle
\enlargethispage{2pt}

\section{Introduction}\label{sec:intro}
\subsection{Gromov ellipticity of irrational Fano varieties}
We work over $\CC$. A smooth algebraic variety $Y$ is \emph{algebraically elliptic} in Gromov's sense if it admits an algebraic vector bundle $p:E\to Y$ and a morphism $s:E\to Y$ such that $s(0_y)=y$ and
\[
 d(s|_{E_y})_{0_y}:T_{0_y}E_y\longrightarrow T_yY
\]
is surjective for every $y\in Y$ \cite{Gro}. The triple $(E,p,s)$ is called a \emph{dominating spray}. Algebraic ellipticity implies the holomorphic Oka property, including approximation and interpolation; see \cite{Zai}. If $Y$ is algebraically elliptic and complete of dimension $n$, by Forstneri\v c's theorem \cite[Theorem~1.6]{For} there is a morphism $\AA^n\to Y$ whose restriction to some Zariski-open subset is smooth and surjective.

Every algebraically elliptic variety is unirational. Conversely, Banecki's theorem on uniform retract rationality implies that every smooth projective rational variety is algebraically elliptic \cite{Ban}.  For smooth projective surfaces, Castelnuovo's theorem therefore identifies algebraic ellipticity with rationality.

In higher dimensions, nonrational elliptic varieties occur among Fano varieties. For example, every smooth cubic hypersurface of dimension at least two is algebraically elliptic \cite{KZcubic}. Smooth cubic threefolds are nonrational by Clemens--Griffiths \cite{CG}, and very general cubic fourfolds are nonrational by Katzarkov--Kontsevich--Pantev--Yu \cite[Theorem~6.8]{KKPY}.
As another example, every smooth quartic double space of dimension at least two
is algebraically elliptic \cite[Theorem~1.1]{DKZ}. Very general quartic double solids and quartic double fourfolds are stably irrational by Voisin \cite{Voi} and Hassett--Pirutka--Tschinkel \cite{HPT}, respectively. For singular examples, Ciliberto and the third-named author constructed rational nodal quartic double solids with uniformly rational, hence algebraically elliptic, small resolutions \cite{CZ}.

Iskovskikh and Manin proved in 1971 that smooth quartic threefolds are nonrational \cite{IM}. Segre's unirational examples are consequently counterexamples to the three-dimensional L\"uroth problem.

The threefolds considered here belong to the family of quartics containing a smooth del Pezzo surface of degree four. Such a surface is an intersection of two quadrics, and a quartic containing it has an equation
\[
 P_1Q_2-P_2Q_1=0,
\]
where the four entries are quadratic forms. For general entries, the quartic has sixteen ordinary double points, precisely the common zeros of the four quadratic forms. Cheltsov proved that a sufficiently general member of this family is nonrational \cite[Theorem~11]{Cheltsov}. His proof uses the degree-four del Pezzo fibration on a small resolution. These threefolds are nevertheless unirational \cite[Remark~12]{Cheltsov}. The total spaces of the associated del Pezzo fibrations are birational to conic bundles over surfaces, see \cite[Corollary~22]{Cheltsov}. These Mori fiber structures are not square birational, since their bases have dimensions one and two. Thus the del Pezzo fibrations are not birationally rigid.

We study diagonal quadratic forms vanishing at the sixteen points $[\pm1:\pm1:\pm1:\pm1:1]$. The right and left kernel constructions define two projective small resolutions, each carrying a degree-four del Pezzo fibration. In Corollary~\ref{cor:nonrationality} we verify the hypotheses of Alexeev's criterion \cite[Theorem~2]{Alexeev} for the diagonal fibration by adapting Cheltsov's argument \cite[\S\,3, proof of Theorem~11]{Cheltsov}. We prove that, on a nonempty Zariski-open subset of this subfamily, both small resolutions and the blow-up of the sixteen nodes admit dominating algebraic sprays. The small resolutions and the blow-ups of the sixteen nodes form two new families of smooth projective algebraically elliptic irrational threefolds; see Theorem~\ref{thm:main} and Remark~\ref{rem:identifying-families} below.

\subsection{The family of quartics and the main theorem}

Our parameter space is
\[
 \mathcal A=\{C\in\operatorname{Mat}_{4\times5}(\CC):C(1,1,1,1,1)^T=0\}
 \simeq\AA^{16}.
\]
For each $C\in\mathcal A$, we define four quadratic forms by
\begin{equation}\label{eq:matrix}
 \begin{pmatrix}P_1\\Q_1\\P_2\\Q_2\end{pmatrix}
 =C\begin{pmatrix}x_0^2\\x_1^2\\x_2^2\\x_3^2\\x_4^2\end{pmatrix},
\end{equation}
and write
\begin{equation}\label{eq:quartic}
 M_C(x)=\begin{pmatrix}P_1(x)&Q_1(x)\\P_2(x)&Q_2(x)\end{pmatrix},
 \qquad f_C(x)=\det M_C(x),\qquad X_C=\{f_C=0\}\subset\PP^4.
\end{equation}
We retain the subscript in $f_C$ even when $C$ is fixed, to distinguish this polynomial from the fibrations denoted by $f$ below.
Since each row of $C$ sums to zero, all four quadrics vanish at the sixteen points
\begin{equation}\label{eq:Sigma}
 \Sigma=\{[\epsilon_0:\epsilon_1:\epsilon_2:\epsilon_3:1]:
                  \epsilon_i\in\{1,-1\}\}
\end{equation}
for every $C\in\mathcal A$.

\begin{theorem}\label{thm:main}
There is a Zariski-open dense subset $\mathcal U_0\subset\mathcal A$ such that for every $C\in\mathcal U_0$ the following hold.
\begin{enumerate}
\item The singular locus of $X_C$ is $\Sigma$, and all sixteen points of $\Sigma$ are ordinary double points.
\item The quartic $X_C$ is non-rational and admits two distinguished small resolutions $X_{R,C}$ and $X_{L,C}$, defined in Section~\ref{sec:genericity}. These resolutions and the common blow-up
\[
 Z_C=\Bl_{\Sigma}X_C
\]
are smooth projective irrational threefolds.
\item The threefolds $Z_C$, $X_{R,C}$, and $X_{L,C}$ are algebraically elliptic.
\end{enumerate}
\end{theorem}

Corollary~\ref{cor:singular-locus} proves assertion (1). Assertion (2) follows from Proposition~\ref{prop:geometry} and Corollary~\ref{cor:nonrationality}. For assertion (3), we treat the common blow-up in Proposition~\ref{prop:Zelliptic} and the two small resolutions in Theorem~\ref{thm:smallres}. Both del Pezzo pencils define morphisms on $Z_C$. On each small resolution only one pencil defines a morphism. On the right small resolution, we first move a point of an exceptional curve by a birational involution. At its image we construct two sprays with independent velocities tangent to the right fiber and a third spray tangent to the left fiber whose velocity lies outside the right vertical tangent plane. We then transport all three sprays back. The left small resolution is treated by symmetry.

\begin{remark}\label{rem:identifying-families}
Interchanging $Q_1$ and $P_2$ transposes $M_C$ and exchanges the right and left constructions without changing $X_C$. More precisely, if $\tau C$ is obtained by exchanging the second and third rows of $C$, then
\[
 X_{R,\tau C}\simeq X_{L,C},\qquad X_{L,\tau C}\simeq X_{R,C}.
\]
These identifications respect the projections to $X_C$ and to $\PP^1$. Thus transposition exchanges the right and left constructions over the full coefficient space $\mathcal A$. For a fixed $C\in\mathcal U_0$, however, the two resolutions are not isomorphic over $X_C$: they correspond to the two local small resolutions at every node and are related by flopping the sixteen exceptional curves.
\end{remark}

\noindent {\bf AI disclosure.}
{Motivated by their work on algebraic ellipticity of cubic threefolds \cite{KZcubic} and quartic double spaces \cite{DKZ}, the authors directed GPT-6 Astra to search for further examples of nonrational algebraically elliptic varieties. This search led to the class considered in this paper. The resulting arguments were revised, expanded where necessary, edited, and streamlined by the authors with assistance from GPT-6 Astra. The authors take responsibility for the mathematical content.}

\section{Restricting the family of quartics}\label{sec:genericity}
\subsection{Del Pezzo fibrations and nonrationality}
For $C\in\mathcal A$, we consider the two closed subschemes
\begin{align}
 X_{R,C}&=\{(x,[u:v]):uP_1+vQ_1=uP_2+vQ_2=0\}
                 \subset\PP^4\times\PP^1,\label{eq:XR}\\
 X_{L,C}&=\{(x,[s:t]):sP_1+tP_2=sQ_1+tQ_2=0\}
                 \subset\PP^4\times\PP^1.\label{eq:XL}
\end{align}
We denote their projections to $X_C$ by $\pi_{R,C},\pi_{L,C}$ and their projections to $\PP^1$ by $\varphi_{R,C},\varphi_{L,C}$. For $C\in\mathcal U$, with $\mathcal U$ defined in Notation~\ref{nota-U} below, Proposition~\ref{prop:geometry} proves that $X_{R,C}$ and $X_{L,C}$ are smooth projective varieties, their projections to $X_C$ are small resolutions, and their fibers over $\PP^1$ are degree-four del Pezzo surfaces. The resolutions are weak Fano, but not Fano, see Remark~\ref{rem:weak-fano}.

We suppress the subscript $C$ when $C$ is fixed.
The fiber of $\varphi_R$ over a point $b=[u:v]\in\PP^1$ is an intersection of two quadrics. For $C\in\mathcal U$, every such surface $\varphi_R^{-1}(b)\subset \PP^4$ contains the sixteen points of $\Sigma$. In the total space $X_{R} \subset\PP^4\times\PP^1$ of the family, these fixed points define the 16 sections $\{p\}\times\PP^1$ of the projection $\varphi_{R}\colon X_{R}\to\PP^1$, which we call the \emph{marked sections}. Each marked section meets the fiber $\varphi_R^{-1}(b)$ in the single point $(p,b)$; these sixteen intersection points are the \emph{marked points} of that fiber.

At a point $x$ of $X_C$ where $M_C(x)\ne0$, the matrix $M_C(x)$ has rank one, so its right kernel determines a unique point of $\PP^1$. Over the open set $V_C=\{x\in X_C:M_C(x)\ne0\}$, the projection $\pi_{R,C}:X_{R,C}\to X_C$, given by $(x,[u:v])\mapsto x$, restricts to an isomorphism
\[
 \pi_{R,C}^{-1}(V_C)\xrightarrow{\ \sim\ }V_C,
\]
whose inverse is
\[
 x\longmapsto(x,[u:v]),\qquad M_C(x)\binom{u}{v}=0,
\]
where $[u:v]\in\PP^1$ is uniquely determined and varies regularly with $x\in V_C$. At a point $p\in\Sigma$, the matrix $M_C(p)$ is zero, so
\[
 \pi_{R,C}^{-1}(p)=\{p\}\times\mathbf P^1.
\]
This curve is a fiber of $\pi_{R,C}$ and a marked section of $\varphi_{R,C}$: the latter projection maps it isomorphically onto $\PP^1$. Under the rank condition below, $\Sigma$ is the whole zero locus of $M_C$, and the sixteen marked sections are the exceptional curves of $\pi_{R,C}$.

\begin{notation}
We write $C=(c_{kl})_{0\le k\le3,\,0\le l\le4}$ and arrange its $i$th column as the $2\times2$ matrix
\[
 \mathsf C_i=\begin{pmatrix}c_{0i}&c_{1i}\\ c_{2i}&c_{3i}\end{pmatrix},
 \qquad \sum_{i=0}^4\mathsf C_i=0.
\]
The coefficient forms of the right pencil are
\[
 a_i^R=uc_{0i}+vc_{1i},\qquad b_i^R=uc_{2i}+vc_{3i},
\]
and those of the left pencil are
\[
 a_i^L=uc_{0i}+vc_{2i},\qquad b_i^L=uc_{1i}+vc_{3i}.
\]
For their minors we write
\[
 D_{ij}^\epsilon=a_i^\epsilon b_j^\epsilon-a_j^\epsilon b_i^\epsilon,
 \qquad \epsilon=R,L.
\]
Let $N_C$ be the unique symmetric $5\times5$ matrix such that
\[
 f_C(x)=z^TN_Cz,\qquad z=(x_0^2,\ldots,x_4^2)^T.
\]
\end{notation}

We introduce next the open subset $\mathcal U\subset\mathcal A$ mentioned before.
The open subset $\mathcal U_0$ in Theorem~\ref{thm:main} will be obtained by shrinking $\mathcal U$ further.

\begin{notation}\label{nota-U}
{\rm Let $\mathcal U\subset\mathcal A$ be the locus on which the following conditions hold:}
\begin{enumerate}[label=(\roman*)]
\item $\operatorname{rank}C=4$ and $\det\mathsf C_i\ne0$ for $0\le i\le4$;
\item for each pencil the ten binary quadratics $D_{ij}^\epsilon$ are nonzero and squarefree, and no two have a common projective zero;
\item $c_{0i}c_{1j}-c_{1i}c_{0j}\ne0$ for all $i<j$; equivalently, the five coefficient forms $uc_{0i}+vc_{1i}$ have distinct projective zeros;
\item every principal minor of $N_C$ of size one, two, or three is nonzero, and each restriction $f_C(e_i+t e_j)$ is a squarefree quartic with nonzero leading and constant terms;
\item the fibers at $[1:1]$ of both pencils are smooth and contain no projective line through a point of $\Sigma$. Here and below, a \emph{line on a del Pezzo surface of degree four} means a projective line in its anticanonical embedding in $\PP^4$.
\end{enumerate}
\end{notation}

These are Zariski-open conditions. Conditions (i)--(iv) are given by the nonvanishing of finitely many minors, discriminants, resultants, leading coefficients, and constant coefficients. For (v), the lines in the fibers at $[1:1]$ form a relative Fano scheme, closed in the product of the coefficient space with $\operatorname{Gr}(1,4)$. Incidence with a fixed point of $\Sigma$ is closed, and its projection to the coefficient space is closed because the Grassmannian is proper. Thus absence of such a line is open. The sign-change group acts transitively on $\Sigma$, so one point of $\Sigma$ suffices.

To show that $\mathcal U$ is Zariski dense in $\mathcal A$, we let
\begin{equation}\label{eq:witnessC}
 C_*=\coeffmatrix.
\end{equation}
In Appendix~\ref{app:geometry}, we verify that $C_*\in\mathcal U$.

\begin{proposition}\label{prop:geometry}
For every $C\in\mathcal U$, the varieties $X_{R,C}$ and $X_{L,C}$ are smooth projective small resolutions of $X_C$. Every fiber of either $\varphi_R$ or $\varphi_L$ is an integral normal degree-four del Pezzo surface. Each of the fibrations $\varphi_R$ and $\varphi_L$ has exactly twenty singular fibers, each a del Pezzo surface of degree four with exactly two $A_1$ singularities; all other fibers are smooth. The sixteen marked sections are pairwise disjoint and meet every fiber at smooth points of that fiber.
\end{proposition}
\begin{proof}
Since $\rk C=4$ and $C(1,1,1,1,1)^T=0$, we have
$\ker C=\CC(1,1,1,1,1)$, so the four quadrics in \eqref{eq:matrix} have common zero locus $\Sigma$. Their differentials are independent there. Hence the four entries of $M_C$ form analytic coordinates near each point of $\Sigma$, and in those coordinates $X_C$ has equation $uz-vw=0$. Thus every point of $\Sigma$ is an ordinary double point of $X_C$.

We analyze the right pencil; the same argument applies to the left pencil after transposing $M_C$, since conditions (i) and (ii) hold for both. For a fixed $\xi=[u:v]\in\PP^1$, we write $a_k=a_k^R(u,v)$ and $b_k=b_k^R(u,v)$. The quartic surface $S_\xi=\varphi_R^{-1}(\xi)$ is defined in $\PP^4$ by
\begin{equation}\label{eq:fiber-quadrics}
 \sum_{k=0}^4a_kx_k^2=\sum_{k=0}^4b_kx_k^2=0.
\end{equation}
By the Jacobian criterion, a point $x\in S_\xi$ is singular if and only if there exists $(\lambda,\mu)\ne(0,0)$ such that
\begin{equation}\label{eq:jacobian-fiber}
 (\lambda a_k+\mu b_k)x_k=0\qquad(0\le k\le4).
\end{equation}
Condition (i) implies that every coefficient column $(a_k,b_k)^T$ is nonzero. A point of $S_\xi$ cannot have only one nonzero coordinate, since the two defining equations would then make the corresponding coefficient column zero. At a singular point, therefore, at least two coefficient columns lie in the one-dimensional kernel of $(\lambda,\mu)$ and are proportional. If the surface $S_\xi$ is singular, there is consequently a pair of indices $i<j$ for which $D_{ij}^R(u,v)=0$. This pair is unique by condition (ii), since distinct minors have no common zero. Every other pair of coefficient columns is linearly independent. Equations~\eqref{eq:jacobian-fiber} then imply that every singular point has $x_k=0$ for $k\notin\{i,j\}$.
By condition (i), neither of the two proportional coefficient columns is zero. Hence, on the coordinate line $l_{i,j}\simeq \PP^1_{\CC}$ spanned by $e_i,e_j$, the two quadric equations defining the del Pezzo surface reduce to the same binary quadratic form
\[
 \alpha x_i^2+\beta x_j^2=0,
 \qquad \alpha\beta\ne0.
\]
On $\PP^1_{\CC}$, this quadratic form has two distinct zeros, which are the only singular points of the fiber.

At either singular point of $S_\xi$, both $x_i$ and $x_j$ are nonzero. We work in the affine chart $x_i=1$. After interchanging the two equations in \eqref{eq:fiber-quadrics} if necessary, we may assume $a_i\ne0$; the proportionality of the nonzero columns implies $a_j\ne0$. The equation $\sum_k a_kx_k^2=0$ can be solved analytically for $x_j$ near the singular point, since its derivative $2a_jx_j$ is nonzero there. The three coordinates $x_k$, $k\notin\{i,j\}$, are therefore local coordinates on this smooth hypersurface, centered at the singular point. Subtracting $b_i/a_i$ times the first equation in \eqref{eq:fiber-quadrics} from the second, we find
\[
 \sum_{k\notin\{i,j\}}\left(b_k-\frac{b_i}{a_i}a_k\right)x_k^2=0.
\]
Every coefficient is nonzero, since the columns indexed by $i$ and $k$ are independent. Thus the surface has a nondegenerate quadratic equation in these three local coordinates, and each singularity is of type $A_1$. To examine the total space, we take a local coordinate $t$ on $\PP^1$ centered at $\xi$. In a local trivialization the coefficients $a_k,b_k$ are functions of $t$. Near either singular point of $S_\xi$, we can still solve the first equation in \eqref{eq:fiber-quadrics} for $x_j$. After this substitution the second equation has the form $F(z,t)=0$, where $z=(x_k)_{k\notin\{i,j\}}$. Since $x_i=1$, at $z=0$ we have
\begin{equation}\label{eq:total-space-smoothing-derivative}
 F(0,t)=b_i(t)-\frac{b_j(t)a_i(t)}{a_j(t)}
       =-\frac{D_{ij}^R(t)}{a_j(t)},\qquad
 \frac{\partial F}{\partial t}(0,0)\ne0.
\end{equation}
The last inequality follows because $a_j(0)\ne0$ and $D_{ij}^R$ has a simple zero at $\xi$, by condition (ii). Hence $X_{R,C}$ is smooth at the point under consideration.

At a smooth point of a fiber, the Jacobian matrix of the two defining equations already has rank two in local coordinates on $\PP^4$. Adding the coordinate $t$ does not decrease this rank, so $X_{R,C}$ is smooth there as well. This applies, in particular, at every point of the marked sections: all its homogeneous coordinates are nonzero, whereas a singular point of a fiber has only two nonzero coordinates. Thus $X_{R,C}$ is smooth everywhere.

The projection $\pi_{R,C}\colon X_{R,C}\to X_C$ is an isomorphism over $X_C\setminus\Sigma$. Since $X_{R,C}$ is smooth, $X_C$ has no singularities outside $\Sigma$. Each fiber of $\varphi_{R,C}$ is a connected complete intersection of two quadrics in $\PP^4$. It is Cohen--Macaulay, hence satisfies Serre's condition $S_2$; its singular locus is finite, so it is regular in codimension one. The fibers are therefore normal by Serre's criterion. A connected normal scheme is irreducible, hence the fibers are integral. By adjunction, $-K_S=\OO_S(1)$ and $K_S^2=4$.

The common resolution is
\begin{equation}\label{eq:Z}
 Z_C=X_{R,C}\times_{X_C}X_{L,C}=\Bl_\Sigma X_C.
\end{equation}
We denote the two projections by
\begin{equation}\label{eq:beta}
\beta_R:Z_C\to X_{R,C},\qquad \beta_L:Z_C\to X_{L,C}
\end{equation}
and the exceptional divisor of $Z_C\to X_C$ over $p\in\Sigma$ by $E_p$. For this divisor and its normal bundle we have
\begin{equation}\label{eq:Ep}
 E_p\simeq\PP_R^1\times\PP_L^1,
 \qquad N_{E_p/Z_C}\simeq\OO(-1,-1).
\end{equation}
The two morphisms
\begin{equation}\label{eq:g}
 g_R=\varphi_R\circ\beta_R\colon Z_C\to\PP^1_R,
 \qquad g_L=\varphi_L\circ\beta_L\colon Z_C\to\PP^1_L,
\end{equation}
together induce the isomorphism
\begin{equation}\label{eq:phiE}
 (g_R,g_L)|_{E_p}:E_p\xrightarrow{\sim}\PP^1\times\PP^1.
\end{equation}
In the analytic coordinates at $p$, the equation of $X_C$ is $uz-vw=0$. The exceptional divisor of the blow-up is the projectivized tangent cone, the smooth quadric $uz-vw=0$ in $\PP^3$. Under the identification $E_p\simeq\PP_R^1\times\PP_L^1$, the restrictions of $\beta_R$ and $\beta_L$ are the first and second projections, respectively. We obtain $X_R$ by contracting the fibers of the first projection on each $E_p$, and $X_L$ by contracting the fibers of the second. Equivalently, $\beta_R$ is the blow-up of the sixteen marked sections in $X_R$, and $\beta_L$ is the analogous blow-up in $X_L$. A fiber of $g_R$ (resp. $g_L$) is therefore the corresponding del Pezzo surface blown up at its sixteen marked points. Formula~\eqref{eq:phiE} also shows that the two fibrations have independent differentials along every exceptional divisor.
\end{proof}

The following consequence of the proof establishes assertion (1) of Theorem~\ref{thm:main}, since $\mathcal U_0\subset\mathcal U$.
\begin{corollary}\label{cor:singular-locus}
For every $C\in\mathcal U$, the singular locus of $X_C$ consists exactly of the sixteen ordinary double points of $\Sigma$.
\end{corollary}
\begin{proof}
The proof of Proposition~\ref{prop:geometry} identifies the points of $\Sigma$ as ordinary double points and proves smoothness of $X_C\setminus\Sigma$.
\end{proof}

\begin{remark}\label{rem:weak-fano}
Let $Y=\PP^4\times\PP^1$ and $C\in\mathcal U$. Since $X_{R,C}\subset Y$ is a smooth complete intersection of two divisors of bidegree $(2,1)$, its canonical bundle satisfies, by adjunction,
\begin{equation}\label{eq:canonical-small-resolution}
 \omega_{X_{R,C}}
 \simeq\bigl(\omega_Y\otimes\OO_Y(4,2)\bigr)|_{X_{R,C}}
 \simeq\OO_Y(-1,0)|_{X_{R,C}}.
\end{equation}
Consequently,
\[
 \omega_{X_{R,C}}^{-1}\simeq\pi_{R,C}^*\OO_{X_C}(1).
\]
Thus the anticanonical divisor is nef and big, with $(-K_{X_{R,C}})^3=4$. It is not ample, since its degree on each exceptional curve $\{p\}\times\PP^1$ is zero. Hence $X_{R,C}$ is weak Fano, but not Fano. The same conclusion holds for $X_{L,C}$.
\end{remark}

For assertion (2) of Theorem~\ref{thm:main}, it remains to check that the diagonal fibration is standard and to compute the Euler characteristic of its total space.

\begin{corollary}\label{cor:nonrationality}
For every $C\in\mathcal U$, the varieties $X_C$, $X_{R,C}$, $X_{L,C}$, and $Z_C=\Bl_\Sigma X_C$ are nonrational.
\end{corollary}
\begin{proof}
We apply Alexeev's criterion \cite[Theorem~2]{Alexeev} in the form stated by Cheltsov \cite[Theorem~23]{Cheltsov}. In $Y=\PP^4\times\PP^1$, let $H_1\subset Y$ be the hypersurface of bidegree $(2,1)$ defined by $uP_1+vQ_1=0$. Condition (iii) ensures that $H_1$ is smooth. Indeed, the five coefficient forms $uc_{0i}+vc_{1i}$ have distinct projective zeros. At a zero of one coefficient, the corresponding quadric has a unique vertex, where the derivative of the defining equation in the base direction is nonzero because that coefficient has a simple zero. At every other point, a derivative in a fiber direction is nonzero. The divisor $H_1$ on $Y$ is ample, and $X_{R,C}$ is the ample divisor on $H_1$ cut out by $uP_2+vQ_2=0$: their associated line bundles are $\OO_Y(2,1)$ and its restriction to $H_1$. Since $H_1$ and $X_{R,C}$ are smooth, the restriction maps on Picard groups are isomorphisms by the Grothendieck--Lefschetz theorem:
\[
 \Pic(Y)\xrightarrow{\sim}\Pic(H_1)\xrightarrow{\sim}\Pic(X_{R,C}),
\]
see \cite[Exp.~XII]{SGA2}. Consequently, $\Pic(X_{R,C})=\ZZ H\oplus\ZZ L$, where $H=\OO_Y(1,0)|_{X_{R,C}}$ and $L=\OO_Y(0,1)|_{X_{R,C}}$. Together with the normality of the fibers established in Proposition~\ref{prop:geometry}, this makes $\varphi_R$ a standard degree-four del Pezzo fibration. To compute the Euler characteristic of $X_{R,C}$, we use its tangent and normal bundle sequences. Its total Chern class is
\[
 c(T_{X_{R,C}})=\frac{(1+H)^5(1+L)^2}{(1+2H+L)^2}.
\]
Since $L^2=0$, the degree-three term is
\[
 c_3(T_{X_{R,C}})=-2H^3-4H^2L.
\]
We also write $H,L$ for the hyperplane classes on $Y$. The class of $X_{R,C}$ in $Y$ is $(2H+L)^2=4H^2+4HL$. Hence
\[
 \int_{X_{R,C}}H^3=\int_{X_{R,C}}H^2L=4,
\]
and therefore
\[
 e(X_{R,C})=\int_{X_{R,C}}c_3(T_{X_{R,C}})
 =-2\cdot4-4\cdot4=-24\notin\{0,-4,-8\},
\]
as required by the criterion. Hence $X_{R,C}$ and its birational models $X_C$, $X_{L,C}$, and $Z_C$ are nonrational.
\end{proof}

\Needspace{12\baselineskip}
In Lemma \ref{lem:twonodallines} below we describe the lines on the singular del Pezzo surfaces in our fibrations. We first introduce the notation.

\begin{notation}\label{nota:S} Fix $C\in\mathcal U$ and a point $\xi=[u_0:v_0]\in\PP^1$ for which $S_\xi=\varphi_{R,C}^{-1}(\xi)$ is singular. Let $i<j$ be the indices of the proportional coefficient columns. For $c=(c_0,\ldots,c_4)$ we choose one of the two rows
\[
 (a_k^R(u_0,v_0))_{k=0}^4,\qquad (b_k^R(u_0,v_0))_{k=0}^4
\]
with nonzero $i$th entry. Since the $i$th and $j$th columns are nonzero and proportional, $c_j\ne0$ as well.

We denote the other row by $h=(h_0,\ldots,h_4)$. The two defining equations of $S_\xi$ in \eqref{eq:fiber-quadrics}, ordered according to the chosen rows, are
\begin{equation}\label{eq:twonodal-original}
 q_c(x)=\sum_{k=0}^4c_kx_k^2=0,\qquad
 q_h(x)=\sum_{k=0}^4h_kx_k^2=0.
\end{equation}
\end{notation}

Recall that a smooth del Pezzo surface of degree 4 in $\PP^4$ contains exactly 16 projective lines.

\begin{lemma}\label{lem:twonodallines}
The singular del Pezzo surface $S_\xi=\varphi_{R,C}^{-1}(\xi)$ contains exactly eight projective lines. For every smooth point $a=[a_0:\cdots:a_4]\in S_\xi$,
\begin{equation}\label{eq:singularlinecriterion}
 a\text{ lies on a line contained in }S_\xi\quad\Longleftrightarrow\quad
 c_i a_i^2+c_j a_j^2=0
\end{equation}
where $i,j$ and $c=(c_0,\ldots,c_4)$ are as in Notation~\ref{nota:S}.
\end{lemma}
\begin{proof}
Since the columns indexed by $i,j$ are proportional, $h_j=(h_i/c_i)c_j$. Subtracting $h_i/c_i$ times the first equation in \eqref{eq:twonodal-original} from the second therefore removes the terms in $x_i^2$ and $x_j^2$ from the second equation. With $z'=(x_k)_{k\notin\{i,j\}}$, the remaining equation is
\begin{equation}\label{eq:twonodal-subtraction}
 d(z'):=q_h(x)-\frac{h_i}{c_i}q_c(x)
 =\sum_{k\notin\{i,j\}}\left(h_k-\frac{h_i}{c_i}c_k\right)x_k^2=0.
\end{equation}
All three coefficients are nonzero, since each column indexed by $k\notin\{i,j\}$ is independent of the $i$th column. Thus $d$ is a nondegenerate quadratic form. The first equation in \eqref{eq:twonodal-original} can be written as
\begin{equation}\label{eq:twonodal-binary}
 c_ix_i^2+c_jx_j^2=e(z'),\qquad
 e(z'):=-\sum_{k\notin\{i,j\}}c_kx_k^2.
\end{equation}
Since $c_ic_j\ne0$, the binary quadratic form on the left of \eqref{eq:twonodal-binary} factors over $\CC$ into two independent linear forms. We use the two factors as coordinates $u,v$ in place of $x_i,x_j$, so that $c_ix_i^2+c_jx_j^2=uv$. The equations are
\begin{equation}\label{eq:twonodalnormal}
 d(z')=0,\qquad uv=e(z').
\end{equation}

We study the lines on $S_\xi$ by means of the linear projection from $\PP^4$ to $\PP^2$ with center the line $\ell=\{z'=0\}$:
\begin{equation}\label{eq:projection-from-line}
 \operatorname{pr}_\ell:\PP^4\dashrightarrow\PP^2,\qquad
 [u:v:z']\longmapsto[z'].
\end{equation}
Since $d(z')=0$ on $S_\xi$, the restriction of this projection to $S_\xi$ takes values in the smooth conic $D=\{d=0\}\subset\PP^2$. We denote this restriction by
\begin{equation}\label{eq:twonodal-projection}
 \rho:S_\xi\dashrightarrow D,\qquad [u:v:z']\longmapsto[z'].
\end{equation}
The line $\ell$ meets $S_\xi$ only at its two $A_1$ points, so $\rho$ is defined at every smooth point of $S_\xi$. A projective line on $S_\xi$ cannot map to a line contained in the smooth plane conic, so it must lie over a point of $d=0$. Above $[z'_0]\in\{d=0\}$, the closure of the fiber is
\[
 uv=e(z'_0)\lambda^2.
\]
It splits into two lines if and only if $e(z'_0)=0$. We show that $d=e=0$ has exactly four solutions in $\PP^2$, all reduced. For the three remaining indices $\{r,s,t\}=\{0,\ldots,4\}\setminus\{i,j\}$, we write
\[
 \alpha_k=h_k-\frac{h_i}{c_i}c_k,\qquad \beta_k=-c_k
 \quad(k\in\{r,s,t\}).
\]
By \eqref{eq:twonodal-subtraction} and \eqref{eq:twonodal-binary}, the equations are
\[
 \alpha_r x_r^2+\alpha_s x_s^2+\alpha_t x_t^2=0,\qquad
 \beta_r x_r^2+\beta_s x_s^2+\beta_t x_t^2=0.
\]
In the squared coordinates $w_k=x_k^2$, $k\in\{r,s,t\}$, these equations form the linear system
\begin{equation}\label{eq:lines-squared-coordinates}
 \begin{pmatrix}\alpha_r&\alpha_s&\alpha_t\\
                 \beta_r&\beta_s&\beta_t\end{pmatrix}
 \begin{pmatrix}w_r\\w_s\\w_t\end{pmatrix}=0.
\end{equation}
Every $2\times2$ minor of this matrix is nonzero: for distinct $k,l\in\{r,s,t\}$,
\[
 \alpha_k\beta_l-\alpha_l\beta_k=c_kh_l-c_lh_k\ne0,
\]
since the corresponding columns of the original coefficient matrix are independent. Hence the matrix has rank two, and its kernel is spanned by a vector $(w_r,w_s,w_t)$ with all three entries nonzero. We normalize this vector by $w_r=1$ and choose square roots $\eta_s^2=w_s$, $\eta_t^2=w_t$. The solutions of $d=e=0$ are exactly
\begin{equation}\label{eq:four-line-parameters}
 [x_r:x_s:x_t]=[1:\varepsilon_s\eta_s:\varepsilon_t\eta_t],
 \qquad \varepsilon_s,\varepsilon_t\in\{1,-1\}.
\end{equation}
They are four distinct points. Each is reduced: in the chart $x_r=1$, the determinant of the Jacobian of the two equations with respect to $x_s,x_t$ is
\[
 4x_sx_t(\alpha_s\beta_t-\alpha_t\beta_s)\ne0.
\] The fiber over each of these four points consists of two lines, and these are all the lines on $S_\xi$.

For a smooth point $a\in S_\xi$, the closure of the fiber over $\rho(a)$ is the union of two lines if and only if $e(z'(a))=0$. Here $z'(a)$ denotes any homogeneous representative of $\rho(a)$; the vanishing condition is independent of this choice. For a homogeneous representative $(a_0,\ldots,a_4)$ of $a$, equation~\eqref{eq:twonodal-binary} reads
\[
 e\bigl((a_k)_{k\notin\{i,j\}}\bigr)=c_i a_i^2+c_j a_j^2.
\]
Thus $e(z'(a))=0$ is equivalent to the condition in \eqref{eq:singularlinecriterion}.
\end{proof}

\begin{remark}\label{rem:lines-specialization}
When a smooth fiber degenerates to $S_\xi$, its sixteen lines coalesce in pairs into the eight lines described above; individual lines can be followed after a finite base change. Four of the eight lines on $S_\xi$ pass through each node, and none passes through both. For each of the four points $q\in\{d=e=0\}$, the conic over $q$ splits into two lines, one through each node, meeting at the smooth point $[0:0:q]$. These are the only intersections between lines through different nodes.
\end{remark}

\subsection{Moving the critical points of \texorpdfstring{$(g_R,g_L)$}{(gR,gL)}}\label{sec:critical-coefficients}
For $C\in\mathcal U$, the resolutions are smooth and their del Pezzo fibers have the geometry described above. We now impose one further condition, needed to construct dominating sprays at every point of $Z_C$. The coefficients satisfying this additional condition form the open subset $\mathcal U_0$ in Theorem~\ref{thm:main}.

Recall that $g_R=\varphi_R\beta_R$ and $g_L=\varphi_L\beta_L$
are morphisms from $Z_C=\Bl_\Sigma X_C$ to $\PP^1$, see \eqref{eq:beta} and \eqref{eq:g}.
We study the critical locus of their product
\[
 \Phi_C=(g_R,g_L):Z_C\longrightarrow\PP^1\times\PP^1.
\]
If $d\Phi_C$ has rank two at $z$, the two tangent planes to the fibers of $g_R$ and $g_L$ through $z$ are distinct and together span $T_zZ_C$. Local sprays with velocities spanning each of these planes would therefore form a dominating family. At a point where the rank of $d\Phi_C$ is smaller, we will first use a birational involution of $Z_C$ to move the point to one where the rank is two. The construction of sprays and their transfer by these involutions will be carried out in Sections~\ref{sec:sprays} and~\ref{sec:critical}; here we establish the required geometric conditions.

\begin{notation}\label{nota:critical}
For $C\in\mathcal U$, let $R_{1,C}$ and $R_{0,C}$ be the sets of points of $Z_C$ at which $d\Phi_C$ has rank one and zero, respectively. We write
\begin{equation}\label{eq:GC}
 R_C=R_{1,C}\cup R_{0,C},\qquad \Omega_C=Z_C\setminus R_C.
\end{equation}
We call $R_C$, equipped with its reduced structure, the \emph{critical locus} of $\Phi_C$, and we call its points \emph{critical points}.

By \eqref{eq:phiE}, the restriction of $\Phi_C$ to each exceptional divisor is an isomorphism onto $\PP^1\times\PP^1$. Hence $R_C$ does not meet these divisors. Using the resolution map, we also regard $R_C$ as a subset of $X_C\setminus\Sigma$.
\end{notation}

Each node $p\in\Sigma$ of $X_C$ determines an involution as follows. A general line through $p$ meets the quartic $X_C$ with multiplicity two at $p$ and in two further points. Exchanging these two points defines a birational involution $\tau_{p,C}$ of $X_C$.
We use the representative of $p$ whose last coordinate is one. Along the line $l$ through $p$ and $x$, the defining quartic form $f_C$ of $X_C$ (see \eqref{eq:quartic}) expands as
\begin{equation}\label{eq:ABuniv}
 f_C(sp+tx)=s^2t^2A_{p,C}(x)+st^3B_{p,C}(x)+t^4f_C(x).
\end{equation}
This defines a quadratic form $A_{p,C}$ and a cubic form $B_{p,C}$. The terms in $s^4$ and $s^3t$ vanish because $p$ is a singular point of $X_C$. For $x\in X_C$, the other intersection point of $l\cap X_C$ is
\begin{equation}\label{eq:tauuniv}
 \tau_{p,C}(x)=[A_{p,C}(x)x-B_{p,C}(x)p]
\end{equation}
where the vector on the right is nonzero.

\begin{definition}\label{def:suitable-node}
Let $x\in R_C$. We say that a node $p\in\Sigma$ is \emph{suitable for $x$} if:
\begin{enumerate}[label=(\alph*)]
\item $A_{p,C}(x)\ne0$;
\item the point $y=\tau_{p,C}(x)$ has all five homogeneous coordinates nonzero and $M_C(y)\ne0$;
\item neither $p$ nor $y$ lies on a line contained in either of the two del Pezzo surfaces through $y$.
\end{enumerate}
In (c), the two surfaces are viewed in $\PP^4$: they are the fibers of del Pezzo fibrations corresponding to the unique lifts of $y$ to $X_{R,C}$ and $X_{L,C}$. These lifts exist because $M_C(y)\ne0$. Both surfaces contain $p$, since every member of either pencil contains $\Sigma$.
\end{definition}

Condition (a) will ensure that the involution is an isomorphism between neighborhoods of $x$ and $y$. Condition (b) will place the lift of $y$ in $\Omega_C$, see \eqref{eq:GC}. Condition (c) concerns the rational curves used later to construct sprays: at $y$ we need curves on each of the two del Pezzo surfaces that do not meet the indeterminacy set of the inverse involution. The role of the two exclusions mentioned in (c) is explained in Section~\ref{sec:critical}.

\begin{proposition}\label{prop:genericcoverage}
There is a nonempty Zariski-open subset $\mathcal U_0\subset\mathcal U$, containing $C_*$ (see \eqref{eq:witnessC}), such that for every $C\in\mathcal U_0$ and every $x\in R_C$, at least one node $p\in\Sigma$ is suitable for $x$ in the sense of Definition~\ref{def:suitable-node}.
\end{proposition}

To prove the proposition, we first describe the critical locus $R_C$ and show that the conditions in Definition~\ref{def:suitable-node} are open, see Lemmas \ref{lem:critical-coordinates}--\ref{lem:suitable-open} below.

\begin{lemma}\label{lem:critical-coordinates}
For a point $x\in X_C\setminus\Sigma$, the rank of $d\Phi_C$ at its lift to $Z_C$ is two, one, or zero according as $x$ has at most one, exactly two, or exactly three zero coordinates. No point of $X_C$ has four zero coordinates. Moreover,
\begin{equation}\label{eq:strata}
 R_C=\bigcup_{|I|=3}C_{I,C},\qquad
 C_{I,C}=X_C\cap\{x_j=0:j\notin I\}.
\end{equation}
The ten curves $C_{I,C}$ are smooth plane quartics. The set $R_{0,C}$ consists of forty distinct points, four on each of the ten coordinate lines. The complement of these points in the union \eqref{eq:strata} is $R_{1,C}$.
\end{lemma}
\begin{proof}
For a fixed $C$, we write $X=X_C$. At a point $x\in X\setminus\Sigma$, the matrix $M_C(x)$ has rank one. Its right and left kernels determine the two parameters of the del Pezzo fibrations. A nonzero rank-one $2\times2$ matrix with prescribed right and left kernels is determined up to multiplication by a scalar. Thus an infinitesimal change of $x$ leaves both parameters unchanged exactly when the corresponding change of $M_C(x)$ is a scalar multiple of $M_C(x)$.

On the affine cone we represent $x$ by a nonzero vector in $\CC^5$ and a tangent vector by $v=(v_0,\ldots,v_4)$. The derivative of the four entries of $M_C$ in direction $v$ is $2C\operatorname{diag}(x)v$, where $\operatorname{diag}(x)$ has diagonal entries $x_0,\ldots,x_4$. Since $\ker C=\CC(1,1,1,1,1)$, the condition just described is equivalent to
\begin{equation}\label{eq:kernelvariation}
 \operatorname{diag}(x)v
 =\lambda(x_0^2,\ldots,x_4^2)^T+\mu(1,1,1,1,1)^T
\end{equation}
for some $\lambda,\mu\in\CC$. Such vectors are tangent to the cone: changing a rank-one matrix by a scalar multiple preserves its determinant to first order.

If all $x_i$ are nonzero, the solutions of \eqref{eq:kernelvariation} are
$v_i=\lambda x_i+\mu/x_i$, a two-dimensional space of solutions. The two directions are independent, since otherwise all $x_i^2$ would be equal and $x$ would belong to $\Sigma$. If exactly one coordinate, say $x_i$, is zero, the $i$th equation forces $\mu=0$. The other components are $v_j=\lambda x_j$, while $v_i$ is arbitrary; the solution space again has dimension two. More generally, if exactly $k\ge1$ coordinates vanish, the solution space has dimension $k+1$.

Vectors differing by a multiple of $x$ represent the same projective tangent vector: this multiple only changes the representative of the point $[x]$. We must therefore subtract one from the dimensions just computed. Since $T_xX$ has dimension three, the resulting ranks of $d\Phi_C$ are two, one, and zero as asserted. Condition (iv) excludes the coordinate vertices, so four coordinates cannot vanish at a point of $X$.

The locus with at least two zero coordinates is the union in \eqref{eq:strata}. To check smoothness of the plane quartics, we write $f_C(x)=z^TN_Cz$, with $z_i=x_i^2$. At a singular point of such a curve, with $J$ the set of indices of its nonzero coordinates, the gradient equations would imply
\[
 (N_C)_{J,J}(x_j^2)_{j\in J}=0.
\]
The vector $(x_j^2)_{j\in J}$ is nonzero, whereas condition (iv) makes $(N_C)_{J,J}$ invertible, a contradiction. Finally, condition (iv) says that the restriction of $f_C$ to each coordinate line has four distinct zeros and that neither endpoint is a zero. Their intersections with $X$ therefore consist of forty distinct points. These are exactly the points with three zero coordinates.
\end{proof}

We next justify the local property of the involution used in condition (a).
\begin{lemma}\label{lem:node-local-isomorphism}
If $x\in X_C$ and $A_{p,C}(x)\ne0$, the formula \eqref{eq:tauuniv} defines an isomorphism between neighborhoods of $x$ and $\tau_{p,C}(x)$ in $X_C$.
\end{lemma}
\begin{proof}
The quadratic form $A_{p,C}$ satisfies
\[
 A_{p,C}(z+\lambda p)=A_{p,C}(z)\qquad(\lambda\in\CC).
\]
This follows by comparing the coefficients of $s^2t^2$ in \eqref{eq:ABuniv} after writing $sp+t(z+\lambda p)=(s+\lambda t)p+tz$. Thus adding a multiple of $p$ does not change the value of the quadratic form $A_{p,C}$. For the representative $y=A_{p,C}(x)x-B_{p,C}(x)p$, homogeneity implies
\begin{equation}\label{eq:Aidentity}
 A_{p,C}(y)=A_{p,C}(x)^3.
\end{equation}
In particular, if $A_{p,C}(x)\ne0$, then $y\ne0$ and $A_{p,C}(y)\ne0$. The formula for $\tau_{p,C}$ is therefore regular both at $x$ and at $y$. The construction by exchanging the two intersections on a line shows that $\tau_{p,C}^2$ is the identity as a rational map. Hence its restrictions near these points are inverse morphisms.
\end{proof}

It remains to show that the existence of a suitable node at every critical point persists when the coefficients vary. For this purpose we define
\[
 \mathscr R=\{(C,x):C\in\mathcal U,\ x\in R_C\}
 \subset\mathcal U\times\PP^4,
 \qquad \pi:\mathscr R\longrightarrow\mathcal U.
\]
By \eqref{eq:strata}, this is a closed subset of $\mathcal U\times\PP^4$. The projection $\pi$ is therefore projective. For each of the sixteen fixed nodes $p$ of $X_C$, let
\[
 \mathscr G_p=\{(C,x)\in\mathscr R:p\text{ is suitable for }x\}.
\]

\begin{lemma}\label{lem:suitable-open}
For each $p\in\Sigma$, the subset $\mathscr G_p$ is Zariski open in $\mathscr R$.
\end{lemma}
\begin{proof}
Condition (a) is the nonvanishing of $A_{p,C}(x)$. On this open set, the five coordinates in \eqref{eq:tauuniv} do not vanish simultaneously, by \eqref{eq:Aidentity}. They therefore define a morphism $(C,x)\mapsto y$. Requiring all coordinates of $y$ to be nonzero, and requiring $M_C(y)\ne0$, are open conditions. This proves openness of (a) and (b).

For (c), we use the following elementary consequence of projectivity: in an algebraic family of surfaces in $\PP^4$, the condition that a specified point lies on a line contained in the surface is closed. Indeed, lines in $\PP^4$ are parametrized by the projective Grassmannian $\operatorname{Gr}(1,4)$. The conditions that a line be contained in the surface and contain the specified point are closed. The image of this incidence subset in the parameter space is closed because the projection is proper.

On the open set defined by (a) and (b), the two points of $\PP^1$ determined by the right and left kernels of $M_C(y)$ vary regularly. Consequently the two del Pezzo surfaces through $y$ form algebraic families there. We apply the preceding observation to each family, first with the point $p$ and then with the point $y$. The four conditions that such a line exist are closed; their complements are the locus defined by (c).
\end{proof}

We can now prove the main proposition. Only the final nonemptiness step uses the explicit computations in the appendices.
\begin{proof}[Proof of Proposition~\ref{prop:genericcoverage}]
Let
\[
 \mathscr B=\mathscr R\setminus\bigcup_{p\in\Sigma}\mathscr G_p.
\]
A point $(C,x)$ belongs to $\mathscr B$ precisely when none of the sixteen nodes is suitable for $x$. By Lemma~\ref{lem:suitable-open}, this is a closed subset of $\mathscr R$. Since $\pi$ is projective, the set $\pi(\mathscr B)$ is closed in $\mathcal U$. We therefore take
\[
 \mathcal U_0=\mathcal U\setminus\pi(\mathscr B).
\]
Thus $C\in\mathcal U_0$ means exactly that every critical point of $\Phi_C$ admits a suitable node. This proves the required property and the openness of $\mathcal U_0$.

For nonemptiness, we use the matrix $C_*$ of \eqref{eq:witnessC}. Its membership in $\mathcal U$ is verified in Appendix~\ref{app:geometry}. The further calculations establish the following two statements:
\begin{enumerate}
\item At every point of the ten plane quartics having exactly two zero coordinates, at least one node satisfies all three conditions (a)--(c), see Appendix~\ref{app:rankone}.
\item At each of the forty points having exactly three zero coordinates, all sixteen nodes satisfy (a)--(c), see Appendix~\ref{app:rankzero}.
\end{enumerate}
By Lemma~\ref{lem:critical-coordinates}, these statements cover every point of $R_{C_*}$. Therefore no point of $\mathscr B$ lies over $C_*$, and $C_*\in\mathcal U_0$. This completes the proof.
\end{proof}

{\bf Convention.} From now on, we fix $C\in\mathcal U_0$ and suppress it from the notation. Thus $X=X_C$, $X_R=X_{R,C}$, $X_L=X_{L,C}$, and $Z=Z_C$; the defining polynomial continues to be denoted by $f_C$.

\section{From rational curves to local sprays}\label{sec:sprays}
Our aim is to construct, at each point of $Z$, finitely many local sprays of rank one whose velocities span its tangent space. We shall construct these sprays from families of rational curves in the fibers of the fibrations $g_R$ and $g_L$ on $Z$.

\subsection{Local sprays and ellipticity}
\begin{definition} Let $V$ be a smooth variety.
A \emph{local rank-one spray} at $z$ with target $V$ is a morphism
\[
 s:U\times\AA^1\to V,\qquad s(u,0)=u,
\]
where $U$ is a Zariski neighborhood of $z$ in $V$. Its \emph{velocity at $z$} is $\frac{\partial s(z,\lambda)}{\partial\lambda}|_{(z,0)}$. A  family of local rank-one sprays is \emph{dominating at $z$} if its velocities at $z$ span $T_zV$.
\end{definition}

It suffices to construct these families separately near each point. The base neighborhood is local, but each spray must be defined for every parameter in $\AA^1$, with values in $V$.

\begin{lemma}\label{lem:localglobal}
If every point of a smooth algebraic variety $V$ has a dominating finite family of local rank-one sprays taking values in $V$, then $V$ is algebraically elliptic.
\end{lemma}
\begin{proof} We
fix a point and a finite local dominating family of rank-one sprays there. Next we apply the rank-one extension construction of \cite[Section~1.2, proof of Lemma~1.3 and Lemma~1.5]{KZ} to each member of the family. After shrinking its local domain, a member of the family extends to a rank-one spray on the whole $V$, defined on a not necessarily trivial line bundle. Its velocity at the fixed point is multiplied by a nonzero scalar, so the extended family still spans the tangent space there.

The domination loci are open and cover $V$, so quasicompactness reduces this cover to the loci of finitely many extended families. Together they form a global dominating finite family. The composition theorem for algebraic sprays \cite[Proposition~1.1 and Corollary~1.2]{KZ} combines this family into a single dominating spray.
\end{proof}

\subsection{Transferring sprays by a birational map}
Suppose a birational map takes a point $x$ to a point $y$ where sprays are already available. Even if the map is an isomorphism near these points, its inverse need not be defined along the whole image of a spray starting near $y$. The following formulation makes the required condition explicit.

For a rational map $\theta:V\dashrightarrow V$, let $\operatorname{Indet}(\theta)$ be the complement of its maximal domain of definition, and write $\operatorname{Dom}(\theta)$ for that domain. When $V$ is fibered by surfaces, $\Indet(\theta)\cap S$ means the set where the map on $V$ is undefined, restricted to the surface $S$. A rational map obtained by restricting the formulas to $S$ may extend further on $S$. Such an extension does not suffice for composing sprays that vary in $V$.

\begin{lemma}\label{lem:conjugation}
Let $\theta:V\dashrightarrow V$ be a birational map of a smooth variety, and suppose it is an isomorphism between neighborhoods of $x$ and $y=\theta(x)$. Set $J=\operatorname{Indet}(\theta^{-1})$. If local sprays $s_i:U_y\times\AA^1\to V\setminus J$ have velocities spanning a subspace $W\subset T_yV$, then, after shrinking a neighborhood $U_x$ of $x$, their conjugates are local sprays at $x$ with velocities spanning $d\theta^{-1}_y(W)$. In particular, the conjugates dominate at $x$ whenever the original family dominates at $y$.
\end{lemma}
\begin{proof} We
shrink $U_x$ so that $\theta(U_x)\subset U_y$, and we let
\begin{equation}\label{eq:transport}
 \widetilde s_i(u,\lambda)=\theta^{-1}\bigl(s_i(\theta(u),\lambda)\bigr).
\end{equation}
By our assumption, $ \widetilde s_i$ is a morphism on $U_x\times\AA^1$, and $\widetilde s_i(u,0)=u$. The derivative of the conjugated spray is
\[\frac{\partial \widetilde s_i(x,\lambda)}{\partial\lambda}\vert_{(x,0)}
 =d\theta^{-1}_y\left(\frac{\partial s_i(y,\lambda)}{\partial \lambda}|_{(y, 0)}\right).
\]
\end{proof}

\subsection{Tangent sections on a smooth or two-nodal surface}
We write the constructions in these two subsections for the right del Pezzo fibration. The same constructions and proofs apply to the left fibration.

We consider a normal del Pezzo surface $S\subset\PP^4$ of degree four.
For a smooth point $a$ of $S$ we write $\mathbb T_aS\subset\PP^4$ for its embedded projective tangent plane and $T_aS$ for its tangent vector space. A hyperplane section of $S$ is a complete intersection of two quadrics in $\PP^3$, hence has arithmetic genus one by adjunction. An irreducible section with an ordinary node therefore has rational normalization. We need a curve with an ordinary node at $a$ that is disjoint from a prescribed finite set. The two branch directions are independent and will be the velocities of the sprays.

\begin{lemma}\label{lem:ordinarycurves}
Let $S$ be a del Pezzo surface from a pencil $\varphi_R:X_R\to\PP^1$, see Proposition~\ref{prop:geometry}, and $a\in S$ be a smooth point away from the lines contained in $S$. For every finite set $K\subset S\setminus\{a\}$, there exists an irreducible anticanonical curve $\gamma\subset S$ with an ordinary node at $a$ as its only singularity, disjoint from $K\cup\Sing S$. Its normalization is a smooth rational curve.
\end{lemma}
\begin{proof}
We write $S=\{Q^{(0)}=Q^{(1)}=0\}\subset \PP^4\times\{b\}$, where $b\in\PP^1$ and $Q^{(0)}, Q^{(1)}$ are quadrics in $\PP^4$. If $\mathbb T_aS$ contains yet another point $a'\in S$, then the line $aa'$ tangent to $Q^{(0)}$ and to $Q^{(1)}$ at $a$ lies on both quadrics, and so is contained in $S$ and passes through $a$. The latter contradicts our assumption. Therefore,
\[
 (\mathbb T_aS\cap S)_\red=\{a\}.
\]
Hyperplanes containing $\mathbb T_aS$ form a pencil with no other base point. Each section has multiplicity at least two at $a$.

In affine coordinates centered at $a$, the defining equations have the form
\begin{equation}\label{eq:quadrics}
 Q^{(i)}(a+z)=\ell_i(z)+q_i(z),\qquad i=0,1,
\end{equation}
where $\ell_i$ is linear and $q_i$ is homogeneous quadratic. We have $T_aS=\{\ell_0=\ell_1=0\}$. On this plane we write
\[
 \overline q_i=q_i|_{T_aS},\qquad i=0,1.
\]
These are quadratic forms in two variables, since $\dim T_aS=2$. A hyperplane in the pencil has equation $\lambda\ell_0+\mu\ell_1=0$, and its section of $S$ has tangent cone $\lambda\overline q_0+\mu\overline q_1=0$ at $a$. The forms $\overline q_0,\overline q_1$ have no common zero in $\PP(T_aS)$: otherwise, for a nonzero representative $v$ of such a zero, both equations $Q^{(i)}(a+tv)=t\ell_i(v)+t^2q_i(v)=0$ would hold for all $t$, so the line through $a$ in direction $v$ would lie on $S$.

Let $\widetilde S$ be the smooth surface formed by resolving the nodes of $S$, if present, and blowing up $a$. We denote the pullback of $\OO_S(1)$ by $H$ and the exceptional curve over $a$ by $e_a$. The transformed pencil has class $D=H-2e_a$. On $e_a\simeq\PP(T_aS)$, the strict transforms cut out the pencil generated by $\overline q_0,\overline q_1$. Since these forms have no common zero, the transformed pencil has no base point on $e_a$. It has no base point on an exceptional curve over a node, since the corresponding node of $S$ was not in $\mathbb T_aS$. Thus the pencil defines a morphism to $\PP^1$.

The restriction of this morphism to $e_a$ has degree two, because it is defined by two sections of $\OO_{\PP^1}(2)$ without a common zero. In particular, a general member of the original pencil has multiplicity exactly two at $a$. The resolutions of the nodes of $S$ are crepant, so
\[
 D^2=0,\qquad D\cdot K_{\widetilde S}=-2.
\]
The morphism $\widetilde S\to\PP^1$ defined by the transformed pencil has Stein factorization
\[
 \widetilde S\xrightarrow{\psi_0}B_a\xrightarrow{\psi_1}\PP^1,
\]
where $B_a$ is a smooth projective curve, $\psi_0$ has connected fibers, and $\psi_1$ is finite. Here $\PP^1$ parametrizes the hyperplanes through $\mathbb T_aS$. If $k=\deg\psi_1$ and $D'$ is the class of a general fiber of $\psi_0$, then $D\equiv kD'$, so $(D')^2=0$ and $K_{\widetilde S}\cdot D'=-2/k$. A general fiber of $\psi_0$ is smooth by generic smoothness. By adjunction,
\[
 2p_a(D')-2=-\frac{2}{k}.
\]
Since the left side is an even integer, $k=1$. Hence the general fiber of the original pencil is connected and smooth, and its genus is zero by adjunction. The degree-two map on $e_a$ is separable. Away from its branch values, a fiber meets $e_a$ at two distinct transverse points; after blowing down, their image is an ordinary node at $a$. General fibers of the pencil also avoid the exceptional curves over the nodes of $S$, which lie in finitely many special fibers, so their images are smooth away from $a$.

Each point of $K\cup\Sing S$ excludes at most one pencil member, since it is not a base point. We choose a general member avoiding these finitely many points. Its image $\gamma\subset S$ has all the required properties.
\end{proof}

\subsection{A relative rank-one spray on \texorpdfstring{$X_R$}{X\_R}}
\begin{definition}\label{def:cS} Over a nonempty open subset $B\subset\PP^1$,
we restrict the del Pezzo fibration $\varphi_R:X_R\to\PP^1$ to
$\cS=\varphi_R^{-1}(B)$ and write $f=\varphi_R|_{\cS}:\cS\to B$.
A local rank-one spray on $X_R$ over $B$ will be called
\emph{relative with respect to $f$} if
each of its orbit curves is contained in a fiber of $f$.
\end{definition}

In Lemma \ref{lem:relativeordinary} we construct such a relative spray.
We use the following notation.

\begin{notation}
Let $f:\cS\to B$ be as above.
We fix $b\in B$ and a smooth point $a\in\cS_b$ lying on no line of $\cS_b$.
By Lemma~\ref{lem:ordinarycurves} we may choose an irreducible anticanonical
curve $\gamma$ with an ordinary node at $a$ as its only singularity and disjoint from $\Sing\cS_b$.
We choose one of its two branches at $a$
and denote by $w\in T_a\cS_b$
a nonzero vector tangent to that branch.
\end{notation}

\begin{lemma}\label{lem:relativeordinary}
Let $F\subset\cS$ be a closed subset disjoint from $\gamma$.
There exists a relative rank-one spray
$s:U\times\AA^1\to\cS\setminus F$ on a Zariski open neighborhood
$U\subset\cS$ of $a$, whose vertical velocity at $a$ is $w$.
\end{lemma}
\begin{proof}
\emph{Constructing the family.}
We first write an equation for the chosen hyperplane section in terms of $a$
and the tangent vector $w$ at $a$. We can identify $\cS_b$
with the intersection $\{Q^{(0)}=Q^{(1)}=0\}\subset\PP^4$.
We work in homogeneous coordinates with $a\in\{x_0\ne0\}$.
If $w$ has components $(w_1,\ldots,w_4)$ in the affine coordinates $x_i/x_0$,
we use the same letter for the vector $(0,w_1,\ldots,w_4)$ when evaluating
the homogeneous quadratic forms $Q^{(0)}$ and $Q^{(1)}$.
The required hyperplane has equation
\begin{equation}\label{eq:tangentsection}
 Q^{(1)}(w)\,dQ^{(0)}|_a-Q^{(0)}(w)\,dQ^{(1)}|_a=0.
\end{equation}
Indeed, this hyperplane contains the tangent plane to $\cS_b$ at $a$, and the quadratic term of its section vanishes in the direction $w$. The two coefficients are not both zero, since at a common zero the corresponding line through $a$ would lie on $\cS_b$, as in the proof of Lemma~\ref{lem:ordinarycurves}. The equation therefore determines the unique member with this tangent direction. The nodal curve $\gamma$ has two distinct tangent directions at $a$, so $w$ determines the chosen branch.

We choose $U$ in the open subset where $f$ is smooth. We extend the nonzero vector $w$
to a regular section $w(u)$ of $\ker df$ on $U$, and we shrink $U$
so that this section never vanishes in $U$. Thus $w(u)$ is tangent to the fiber $\cS_{f(u)}$ through $u$.
The pulled-back family of del Pezzo surfaces $\cS\times_B U\to U$
has the diagonal section $u\mapsto(u,u)$.
 For each $u\in U$, we apply equation~\eqref{eq:tangentsection} to the quadrics that define $\cS_{f(u)}$ instead of $\cS_b$, to the point $u$ instead of $a$, and to the vector $w(u)$ instead of $w$. Let $H_u$ be the hyperplane defined by that equation. The curves $\gamma_u=\cS_{f(u)}\cap H_u$ form a family through the moving point $u$, with prescribed tangent direction $w(u)$, and $\gamma_a=\gamma$.

We blow up the diagonal section within  $\cS\times_B U$. Let $H$ be the relative hyperplane class and $E$ be the exceptional divisor. The strict transform $\gamma'$ of $\gamma$ is its smooth rational normalization, as in Lemma~\ref{lem:ordinarycurves}. Since $\gamma$ avoids the singularities of $\cS_b$, the pulled-back family is smooth over $U$ near $\gamma$. After blowing up the diagonal section, the resulting family remains smooth over $U$ near $\gamma'$. The transformed equation is a section of $\OO(H-2E)$; its restriction to $\Bl_a(\cS_b)$ defines $\gamma'$ as an effective Cartier divisor. Let $\cC$ be the zero scheme of the transformed equation. Its fiber over $a$ is $\gamma'$. The local criterion for flatness shows that $\cC\to U$ is flat near $\gamma'$, and hence smooth there, since $\gamma'$ is smooth. The complement in $\cC$ of the locus where this morphism is smooth is closed and disjoint from $\gamma'$. Since $\cC\to U$ is proper, the image of this complement is a closed subset of $U$ not containing $a$. Over its complement we have a smooth proper family of curves $\pi:\cC\to U$ with $\cC_a=\gamma'$. After a further shrinking of $U$, its fibers are connected and have genus zero, as does $\gamma'$.

On the blown-up surface $\Bl_u(\cS_{f(u)})$, the exceptional curve over $u$ is
\[
 E_u=\PP(T_u\cS_{f(u)})\simeq\PP^1;
\]
its points represent tangent directions to $\cS_{f(u)}$ at $u$.
After shrinking $U$, the strict transform $\cC_u$ of $\gamma_u$ meets $E_u$ in two distinct points, corresponding to the two branches of $\gamma_u$ at $u$.
Equation~\eqref{eq:tangentsection} ensures that the point $[w(u)]\in E_u$ belongs to $\cC_u$.
We distinguish this point and set $\sigma(u)=[w(u)]$.
Since $w(u)$ is a regular nowhere-vanishing section of $\ker df$, these points vary regularly with $u$ and define a section $\sigma:U\to\cC$.
At $u=a$, the point $\sigma(a)$ corresponds to the chosen branch of $\gamma$.
 Let $\Delta=\{(u,u):u\in U\}\subset\cS\times_B U$, and denote the blow-down by
\[
 \rho:\Bl_\Delta(\cS\times_B U)\longrightarrow\cS\times_B U.
\]
Writing $\operatorname{pr}_1:\cS\times_B U\to\cS$ and
$\operatorname{pr}_2:\cS\times_B U\to U$ for the two projections, we define
\[
 \nu=\operatorname{pr}_1\circ\rho|_{\cC}:\cC\longrightarrow\cS,
 \qquad \pi=\operatorname{pr}_2\circ\rho|_{\cC}.
\]
Since $\rho(\sigma(u))=(u,u)$, we have $\nu(\sigma(u))=u$.
The restriction of $\nu$ to $\cC_a=\gamma'$ is the normalization map onto $\gamma$.
Finally, the identity $f\circ\operatorname{pr}_1=(f|_U)\circ\operatorname{pr}_2$
on the fiber product implies $f\nu=(f|_U)\pi$.

\emph{Parametrizing the curves.}
We now choose a coordinate on each curve $\cC_u\simeq\PP^1$, varying regularly with $u$, in which the point $\sigma(u)$ has coordinate zero. Each curve $\cC_u$ is isomorphic to $\PP^1$ because it is a smooth proper curve of genus zero with the point $\sigma(u)$, see \cite[Section~53.10, Proposition~53.10.4]{Stacks}.
The image $\sigma(U)$ is an effective Cartier divisor on $\cC$.
Let $L=\OO_{\cC}(\sigma(U))$ be its associated line bundle, and let $s_0$ be the section of $L$ whose zero divisor is $\sigma(U)$.
On each curve $\cC_u$, the restriction $L_u=L|_{\cC_u}$ has degree one, so
\[
 L_u\simeq\OO_{\PP^1}(1),\qquad
 \dim H^0(\cC_u,L_u)=2,\qquad H^1(\cC_u,L_u)=0.
\]
Thus any basis of $H^0(\cC_u,L_u)$ defines homogeneous coordinates on $\cC_u$.

To choose these sections consistently as $u$ varies, we work with $E_\pi=\pi_*L$.
By the cohomology and base change lemma \cite[Lemma~30.22.1, Tag~07VK]{Stacks}, this is a vector bundle of rank two on $U$, and its fiber at $u$ is $H^0(\cC_u,L_u)$.
The restriction of $s_0$ to every fiber is a nonzero vector in this two-dimensional space.
After shrinking $U$ around $a$, we can complete it to a basis $s_0,s_1$ of $E_\pi$.
Equivalently, $s_0,s_1$ are sections of $L$ whose restrictions form a basis of $H^0(\cC_u,L_u)$ for every $u\in U$.
They have no common zero, since two independent sections of $\OO_{\PP^1}(1)$ have no common zero.
Consequently they define a morphism
\[
 \Phi:\cC\longrightarrow U\times\PP^1,
 \qquad x\longmapsto\bigl(\pi(x),[s_0(x):s_1(x)]\bigr).
\]

On every fiber, $\Phi$ is the isomorphism to $\PP^1$ given by a basis of the sections of a degree-one line bundle.
The morphism $\Phi$ is proper and has finite fibers, hence is finite. Its restriction to the generic fiber over $U$ is an isomorphism, so it is birational. Since $U\times\PP^1$ is smooth, hence normal, $\Phi$ is an isomorphism.
Moreover, $s_0(\sigma(u))=0$ and $s_1(\sigma(u))\ne0$, so $\Phi(\sigma(u))=(u,[0:1])$.

\emph{Keeping the curves disjoint from a closed set.}
Since $\nu(\cC_a)=\gamma$ is disjoint from $F$, the closed subset
$\pi(\nu^{-1}(F))\subset U$ does not contain $a$. Removing it ensures
$\nu(\cC_u)\cap F=\varnothing$ for every $u\in U$.

We shall also use the following consequence of this construction.
Let $M$ be the union of the sixteen marked sections and let
$\beta:\widehat{\cS}=\Bl_M\cS\to\cS$ be the blow-up along $M$.
If $\gamma\cap M=\varnothing$, the same argument lets us shrink $U$
so that $\nu(\cC)\subset\cS\setminus M$. Thus $\nu$ lifts uniquely to
$\widehat\nu:\cC\to\widehat{\cS}$ with $\beta\widehat\nu=\nu$.
For any closed subset $\widehat F\subset\widehat{\cS}$ disjoint from
$\widehat\nu(\cC_a)$, removing the closed subset
$\pi(\widehat\nu^{-1}(\widehat F))$ of $U$ ensures that all the curves
$\widehat\nu(\cC_u)$ are disjoint from $\widehat F$.

\emph{Defining the spray.}
In a local trivialization of the projective-line bundle, we take $\sigma$ as the constant section $0$ and we choose the constant section $\infty$ disjoint from it. With $\iota:U\times\AA^1\hookrightarrow\cC$ the resulting open immersion, the map
\[
 s(u,\lambda)=\nu(\iota(u,\lambda))
\]
is a relative local spray. The central normalization is immersed at $\sigma(a)$, so $\partial_\lambda s(a,0)=cw$ for some $c\in\CC^*$. After replacing $s(u,\lambda)$ by $s(u,\lambda/c)$, the velocity is $w$. When $\gamma\cap M=\varnothing$, the morphism $\widehat\nu\circ\iota$ is similarly a local spray on $\widehat{\cS}$, after identifying $U$ with $\beta^{-1}(U)$.
\end{proof}

\begin{remark}
The construction in Lemma~\ref{lem:relativeordinary} includes a smooth proper family $\pi:\cC\to U$ with evaluation map $\nu:\cC\to\cS$. After shrinking $U$, its complete images avoid any prescribed closed subset disjoint from the central curve $\gamma$. If $\gamma$ avoids the marked sections, the family lifts to their blow-up; after a further shrinking, its complete images avoid any prescribed closed subset disjoint from the lifted central curve.
\end{remark}

\section{Sprays along the fibers of the two fibrations}\label{sec:verticalsprays}
We shall prove that, whenever $g_R$ (resp. $g_L$) in \eqref{eq:g} is smooth at $z\in Z$, there exist two relative local rank-one sprays whose velocities at $z$ span $\ker(dg_R)_z$ (resp. $\ker(dg_L)_z$). At a point $z\in\Omega$, see \eqref{eq:GC}, these four rank-one sprays form a dominating local family on $Z$, because the two kernels at $z$ are distinct planes in $T_zZ$.

For the right fibration, the construction in Section~\ref{sec:sprays} applies directly when $\beta_R(z)$ lies on no line and is not a marked point of the del Pezzo surface $\varphi_R^{-1}(g_R(z))$. The analogous statement holds for the left fibration, with $\beta_L(z)$ and $\varphi_L^{-1}(g_L(z))$. In the remaining cases, we work with each fibration separately: we move $z$ within its fiber to a point $z'$ whose image has these properties, and transfer the sprays back. The point $z'$ may be different for the two fibrations.

\subsection{Points away from the lines and the exceptional curves}
\begin{lemma}\label{lem:vertical-away-lines}
Suppose that $g_R$ is smooth at $z\in Z$ and that $\beta_R(z)$ lies on no line and is not a marked point of the del Pezzo surface $\varphi_R^{-1}(g_R(z))$. Then there exist two local rank-one sprays relative to $g_R$ whose velocities at $z$ span $\ker(dg_R)_z$. The analogous assertion holds for the left fibration.
\end{lemma}
\begin{proof}
Lemma~\ref{lem:ordinarycurves} supplies a curve with an ordinary node at $\beta_R(z)$ that avoids the marked points and the singularities of the surface. We apply Lemma~\ref{lem:relativeordinary} separately to its two branches. After shrinking their base neighborhoods, the images of the resulting sprays avoid the marked sections, so the sprays lift to $Z$, as explained in the proof of that lemma. Their velocities are independent because the two branches have distinct tangent directions and $\beta_R$ is an isomorphism near $z$. Both velocities belong to the two-dimensional space $\ker(dg_R)_z$, and hence span it. The same argument applies to the left fibration.
\end{proof}

\subsection{The relocation theorem}
We now consider a point $z\in Z$ at which $g_R$ is smooth, but whose image $\beta_R(z)$, with $\beta_R$ as in \eqref{eq:beta}, lies on a line in the corresponding del Pezzo surface, or which belongs to an exceptional curve of $\beta_R$. To construct two relative local sprays whose velocities span $\ker(dg_R)_z$, we use a birational involution preserving the fibers. The same argument applies to the left fibration. We state the required result in Theorem~\ref{thm:relocation}; its proof is given in Section~\ref{sec:relocation}.

We fix a pencil of del Pezzo surfaces $f:\cS\to B$ as in Definition \ref{def:cS}.
Let $s_i: B\to\cS$, $i=1\ldots,16$ be its marked sections, $\mathcal Z=\bigcup s_i(B)\subset\cS$,
and set
\begin{equation}\label{eq:gZ}
 \beta:Z_B=\Bl_{\mathcal Z}\cS\longrightarrow\cS,
 \qquad g=f\beta: Z_B\longrightarrow B.
\end{equation}
The required birational involution preserves $g$ and has only finitely many points of indeterminacy on the fiber $(Z_B)_b$, where $b=g(z)\in B$.

\begin{theorem}\label{thm:relocation}
Suppose $g$ in \eqref{eq:gZ} is smooth at $z\in(Z_B)_b$. There is a birational involution $\theta:Z_B\dashrightarrow Z_B$ satisfying $g\theta=g$ such that $\theta$ is an isomorphism between neighborhoods of $z$ and $y=\theta(z)$, the point $\beta(y)$ is smooth, lies on no line contained in $\cS_b$, is not a marked point, and
\begin{equation}\label{eq:finiteindet}
 \Indet(\theta^{-1})\cap(Z_B)_b\quad\text{is finite}.
\end{equation}
\end{theorem}

\subsection{Constructing the sprays}
We apply Theorem~\ref{thm:relocation} to move the point $z$, and then choose rational curves on $\cS_b$ whose lifts to $Z_B$ avoid the indeterminacy locus of the inverse involution.

\begin{corollary}\label{cor:verticalsprays}
At every point $z\in Z_B$ where $g$  in \eqref{eq:gZ} is smooth, there are two local $g$-relative rank-one sprays with target $Z_B$ whose velocities span $\ker(dg)_z$. Every point of $\Omega$ has a dominating finite family of local rank-one sprays with target $Z$.
\end{corollary}
\begin{proof}
For $\theta$ as in Theorem~\ref{thm:relocation}, we let $y=\theta(z)$ and $J=\Indet(\theta^{-1})$. The map $\theta$ is a local isomorphism at $z$, and $\beta(y)$ lies off the marked points and lines on the del Pezzo surface, see Theorem~\ref{thm:relocation}. On the surface $S=\cS_b$, the points to avoid are $\beta(J\cap(Z_B)_b)$, the marked points, and the points of $\Sing S$. Their union is finite and does not contain $\beta(y)$: the inverse map is regular at $y$, and $\beta$ is an isomorphism over $\beta(y)$.

We choose a curve with an ordinary node at $\beta(y)$ that avoids this finite set, see Lemma~\ref{lem:ordinarycurves}. Its lift avoids $J$. We apply the construction in Lemma~\ref{lem:relativeordinary} separately to its two branches. After shrinking the base neighborhoods, the complete images of both resulting proper families avoid $J$. Restricting each family to the complement of a section at infinity defines two local sprays with images in $\Dom(\theta^{-1})$. Their velocities at $y$ are independent because the two branches of an ordinary node have distinct tangent directions.

The two velocities span $\ker(dg)_y$. By Lemma~\ref{lem:conjugation} we can transport both sprays to $z$. Since $g\theta=g$, the differential $d\theta_z$ identifies $\ker(dg)_z$ with $\ker(dg)_y$, so the transported velocities span $\ker(dg)_z$.

At a point of $\Omega$, both $g_R$ and $g_L$ are smooth. Their differentials are linearly independent, so
\begin{equation}\label{eq:tangentspan}
 \ker(dg_R)_z+\ker(dg_L)_z=T_zZ.
\end{equation}
We apply the first assertion separately to the two fibrations.
\end{proof}

\section{Node involutions on the critical locus \texorpdfstring{$R_C$}{R\_C}}\label{sec:critical}
At every point of the open dense subset $\Omega$ of $Z$, we now have a dominating family of local sprays.
It remains to treat its complement $R=R_1\cup R_0=Z\setminus\Omega$, see \eqref{eq:GC}.
At these points, the tangent directions along the two fibrations need not span the tangent space of $Z$.

For each such point $x$, there is a suitable node $p$, see Proposition~\ref{prop:genericcoverage}. The lifted involution $\widehat\tau_p$ on $Z$ takes $x$ to a point $y\in\Omega$ and is an isomorphism near $x$ and $y$. We must choose rank-one local sprays at $y$ whose entire images are disjoint from the indeterminacy locus of the inverse involution. We first describe a closed set containing this locus, and then construct curves disjoint from it.

\subsection{Where the node involution is defined}\label{ss:5.1}
As in \eqref{eq:ABuniv}, for a node $p\in\Sigma$ of $X$ we define $A_p,B_p$ by
\begin{equation}\label{eq:AB}
 f_C(sp+tx)=s^2t^2A_p(x)+st^3B_p(x)+t^4f_C(x).
\end{equation}
For a general line through $p$, its intersection with the quartic consists of the double point $p$ and two residual points. The rational involution interchanging the residual points is
\begin{equation}\label{eq:tau}
 \tau_p(x)=[A_p(x)x-B_p(x)p]
\end{equation}
and it induces a rational involution $\widehat\tau_p$ on $Z$. Let $\cL_p$ be the union of lines on $X$ through $p$.

\begin{lemma}\label{lem:linefiniteness}
For the fixed member $C\in\mathcal U_0$, $\cL_p$ is a finite union of lines. On $X\setminus(\Sigma\cup\cL_p)$, the map $\tau_p$ lifts to a morphism with values in $Z$.
\end{lemma}
\begin{proof}
A second point $q$ on a line $\ell\subset X$ through $p$ determines a parametrization of $\ell$. In homogeneous representatives of $p$ and $q$, its points are $[sp+tq]$, where $[s:t]\in\PP^1$. The entries of $M=M_C$, see \eqref{eq:quartic}, are quadratic forms, and $M(p)=0$. Hence there are constant $2\times2$ matrices $A,B$ such that
\[
 M(sp+tq)=stA+t^2B=t(sA+tB).
\]
Dividing the four entries by their common factor $t$ leaves the matrix $sA+tB$, whose entries are linear in $s,t$. Since $\ell\subset X=\{\det M=0\}$, we have $\det(sA+tB)=0$ identically. Thus these matrices have rank at most one. They are not all zero, since the zero locus of $M$ is the finite set $\Sigma=X_{\rm sing}$.

The projective classes of the nonzero matrices $sA+tB$ have closure a line or a point on the quadric of rank-one matrices
\[
 \left\{\left[\begin{pmatrix}a&b\\c&d\end{pmatrix}\right]:ad-bc=0\right\}\subset\PP^3.
\]
This is the Segre quadric $\PP^1\times\PP^1$: a rank-one matrix is a product of a column vector and a row vector. Along one ruling its right kernel is constant, and along the other its left kernel is constant. Every line on the quadric belongs to one of these rulings. If the projective image is a point, both kernels are constant.

Suppose, for example, that the right kernel is the constant point $[r_0:r_1]$. Then
\[
 M(x)\binom{r_0}{r_1}=0\qquad(x\in\ell),
\]
so the line $\{(x,[r_0:r_1]):x\in\ell\}\subset X_R$ lies in the fiber $\varphi_R^{-1}([r_0:r_1])$. If the left kernel is constant, the corresponding lift lies in a fiber of $\varphi_L$ on $X_L$. Thus every line on $X$ through $p$ lifts to a line in a fiber of $\varphi_R$ on $X_R$ or of $\varphi_L$ on $X_L$.

We treat the right pencil first. For $b\in\PP^1$, we identify the fiber of $\varphi_R$ with the corresponding intersection of two quadrics $S_b\subset\PP^4$. The set of points $b$ in the base $\PP^1$ for which $S_b$ contains a line through $p$ is closed. Indeed, the pairs $(b,\ell)$ with $p\in\ell\subset S_b$ form a closed subset of $\PP^1\times\operatorname{Gr}(1,4)$, where $\operatorname{Gr}(1,4)$ is the projective Grassmannian of lines in $\PP^4$. Its projection to the base $\PP^1$ is therefore closed. By condition (v), the fiber over $[1:1]$ contains no line through $p$, so this closed subset is proper and hence finite. The same argument applies to $\varphi_L:X_L\to\PP^1$.

Each smooth del Pezzo surface of degree four contains only finitely many lines, and the same holds for the two-nodal surfaces, see Lemma~\ref{lem:twonodallines}. Every line on $X$ through $p$ lifts to a line in one of the finitely many fibers just identified. Consequently $\cL_p$ is a finite union of lines.

We now show that $\tau_p$ lifts to a morphism from $V=X\setminus(\Sigma\cup\cL_p)$ to $Z$. First, formula~\eqref{eq:tau} defines a morphism on $V$. Indeed, if $A_p(x)x-B_p(x)p=0$ for $x\in V$, then $x\ne p$ implies $A_p(x)=B_p(x)=0$. Together with $f_C(x)=0$, equation~\eqref{eq:AB} would then imply that the whole line through $p$ and $x$ is contained in $X$, contrary to $x\notin\cL_p$.

To lift $\tau_p:V\to X$ to $\Bl_pX$, we need a local calculation only at points $x_*\in V$ for which $\tau_p(x_*)=p$. At such a point, formula~\eqref{eq:tau} implies that $A_p(x_*)=0$ and $B_p(x_*)\ne0$. In an affine chart containing $p$ and $x_*$, we write their homogeneous representatives as
\[
 p=(1,p_1,\ldots,p_4),\qquad x=(1,x_1,\ldots,x_4).
\]
Here $p_j$ denotes a coordinate of the fixed node $p$. In these coordinates, the image $y=\tau_p(x)$ satisfies
\[
 y_j=\frac{A_p(x)x_j-B_p(x)p_j}{A_p(x)-B_p(x)},\qquad
 y_j-p_j=\frac{A_p(x)(x_j-p_j)}{A_p(x)-B_p(x)}
 \quad(1\le j\le4).
\]
The denominator is nonzero at $x_*$, and at least one difference $x_j-p_j$ is nonzero there because $x_*\ne p$. After restricting to a neighborhood of $x_*$, both are invertible regular functions. In the affine coordinates $y_j$, the differences $y_j-p_j$ generate the ideal of $p$, so its pullback is the principal ideal generated by $A_p$. Since $V$ is integral and $A_p$ is not identically zero, this ideal is invertible. The map therefore lifts, by the universal property of the blow-up, to $\Bl_pX$ near $x_*$. Away from the inverse image of $p$, the lift is obtained from the isomorphism $\Bl_pX\setminus E_p\simeq X\setminus\{p\}$. These local lifts agree and define a morphism $V\to\Bl_pX$.

Finally, no point of $V$ is mapped to another node $q\in\Sigma\setminus\{p\}$. Otherwise, for $x\in V$ with $\tau_p(x)=q$, formula~\eqref{eq:tau} shows that $p,q,x$ lie on the same line. These are three distinct points. If that line were not contained in $X$, its intersection multiplicity with $X$ would be at least two at each node and at least one at $x$, so it would be at least
\[
 2+2+1=5>4.
\]
This contradicts the degree of the quartic. Thus the line would be contained in $X$, contrary to $x\notin\cL_p$. The image of the lift consequently avoids the points over all nodes other than $p$. The blow-up at the remaining nodes is an isomorphism over this image, so the lift extends to a morphism $V\to Z=\Bl_\Sigma X$.
\end{proof}

\subsection{Constructing sprays at a critical point}
The conditions imposed on a suitable node in Section~\ref{sec:critical-coefficients} ensure that the curves used for the sprays can be chosen disjoint from $\Sigma\cup\cL_p$.

\begin{notation} We write $\rho:Z\to X$ for the resolution. Each $v\in X\setminus\Sigma$ has a unique lift $\widetilde v$ to $Z$. For a nonnodal point $y$, let $S_R(y),S_L(y)\subset\PP^4$ be the images under the first projection of the fibers of $\varphi_R,\varphi_L$ through its lifts to $X_R,X_L$. Each such projection is an isomorphism onto an intersection of two quadrics; the corresponding fiber on $Z$ is its blow-up at the sixteen marked points.
\end{notation}

\begin{proposition}\label{prop:ordinaryavoid}
Let $p\in\Sigma$ and $x\in X\setminus\Sigma$ satisfy $A_p(x)\ne0$. For $y=\tau_p(x)$, assume that $y\notin\Sigma$ and $\widetilde y\in\Omega$. Suppose that neither $S_R(y)$ nor $S_L(y)$ contains a line through $p$ or through $y$. Then there is a dominating finite family of local rank-one sprays at $\widetilde x$ taking values in $Z$.
\end{proposition}
\begin{proof}
For the homogeneous representative $Y=A_p(x)x-B_p(x)p$, $A_p(Y)=A_p(x)^3\ne0$ by \eqref{eq:Aidentity}. Thus $y$ is well defined, and $\tau_p$ is a local isomorphism at $x$ and $y$. Since both points are nonnodal, $\widehat\tau_p$ is a local isomorphism at their lifts $\widetilde x,\widetilde y$. Moreover, $y\notin\cL_p$, since a line contained in $X$ through $p$ would force $A_p$ to vanish along it.

We work first on $S=S_R(y)$. This is a smooth or two-nodal del Pezzo surface by Proposition~\ref{prop:geometry}, and it is smooth at $y$ because $\widetilde y\in\Omega$. The finite union $\cL_p$ has finite intersection with $S$: an infinite intersection with one of its lines would imply that the line is contained in $S$, contrary to the hypothesis. Thus
\[
 K=(S\cap\cL_p)\cup\Sigma\cup\Sing S
\]
is a finite set disjoint from $y$.

We choose a curve on $S$ with an ordinary node at $y$ that avoids $K$, see Lemma~\ref{lem:ordinarycurves}. Its lift to $Z$ is disjoint from $\rho^{-1}(\Sigma\cup\cL_p)$. On the complement of this closed subset, $\widehat\tau_p$ is a morphism, see Lemma~\ref{lem:linefiniteness}. Since it is an involution, this complement is also contained in the domain of $\widehat\tau_p^{-1}$. We apply the construction in Lemma~\ref{lem:relativeordinary} separately to the two branches of the chosen curve, and we lift the resulting families to $Z$. After shrinking their base neighborhoods, their complete images remain disjoint from $\rho^{-1}(\Sigma\cup\cL_p)$. The resulting two local $g_R$-relative rank-one sprays at $\widetilde y$ have velocities spanning $\ker(dg_R)_{\widetilde y}$ and images lying in the domain of $\widehat\tau_p^{-1}$.

The same construction on $S_L(y)$ defines two local $g_L$-relative sprays. Since $\widetilde y\in\Omega$, the two vertical planes together span $T_{\widetilde y}Z$, by \eqref{eq:tangentspan}. The four sprays conjugated by $\widehat\tau_p$ take values in $Z$ and dominate at $\widetilde x$ by Lemma~\ref{lem:conjugation}.
\end{proof}

\begin{proposition}\label{prop:coverage}
For every $z\in R_1\cup R_0$, the point $x=\rho(z)$ satisfies the hypotheses of Proposition~\ref{prop:ordinaryavoid} for some $p\in\Sigma$.
\end{proposition}
\begin{proof}
Let $z\in R_1\cup R_0$ and set $x=\rho(z)$. Since $C\in\mathcal U_0$, there exists $p\in\Sigma$ such that $(C,x)\in\mathscr G_p$, see Proposition~\ref{prop:genericcoverage}. Conditions (a)--(c) of Definition~\ref{def:suitable-node}, together with Proposition~\ref{prop:geometry} and Lemma~\ref{lem:critical-coordinates}, imply all the hypotheses of Proposition~\ref{prop:ordinaryavoid}.
\end{proof}

\Needspace{10\baselineskip}

\subsection{Ellipticity of the common blow-up}\label{sec:commonblowup}

\begin{proposition}\label{prop:Zelliptic}
The threefold $Z$ is smooth, projective, nonrational, and algebraically elliptic.
\end{proposition}
\begin{proof}
The variety $Z$ is the smooth projective blow-up of the sixteen nodes, as described in Proposition~\ref{prop:geometry} and \eqref{eq:Z}. It is nonrational, see Corollary~\ref{cor:nonrationality}.

To prove ellipticity, we use the decomposition $Z=\Omega\sqcup R_1\sqcup R_0$. The dominating local families of sprays with values in $Z$ exist on $\Omega$ by Corollary~\ref{cor:verticalsprays} and at every point of $R=R_1\cup R_0$ by Propositions~\ref{prop:coverage} and~\ref{prop:ordinaryavoid}. The variety $Z$ consequently admits a global dominating algebraic spray, by Lemma~\ref{lem:localglobal}.
\end{proof}

\section{Moving a point by a reflection on a conic}\label{sec:relocation}
To prove Theorem~\ref{thm:relocation}, we return to the maps $f:\cS\to B$, $\beta:Z_B\to\cS$, and $g=f\beta$ of \eqref{eq:gZ}. We fix $b\in B$ and a point $z\in(Z_B)_b=g^{-1}(b)$ at which $g$ is smooth. We seek a birational involution that moves $z$ to a point lying on no line or exceptional curve, is an isomorphism near $z$, and has only finitely many points of indeterminacy on the fiber $(Z_B)_b$.

The construction uses one del Pezzo fibration at a time, over a neighborhood $B$ of $b=g(z)$. The fiber over $b$ may be singular away from $z$.

\subsection{An auxiliary section of the del Pezzo fibration}\label{sec:auxiliarysection}
\begin{notation}\label{nota:azq}
We write $S=\cS_b$ and $a=\beta(z)$. We choose an auxiliary smooth point $q\in S$, distinct from $a$, such that the hyperplane containing $\mathbb T_qS$ and $a$ cuts a curve $C(q,a)$ with an ordinary node at $q$, smooth at $a$. If $a$ is marked, the point $z$ on the exceptional curve over $a$ represents a tangent direction at $a$. We require the tangent line to $C(q,a)$ at $a$ to differ from this direction; equivalently, its strict transform meets the exceptional curve at a point different from $z$.
\end{notation}

\begin{lemma}\label{lem:movingcenter}
Let $S$ be a surface as in Lemma~\ref{lem:ordinarycurves}, let $a$ be a smooth point of $S$, let $K\subset S\setminus\{a\}$ be finite, and fix a direction $d\in\PP(T_aS)$. Then there exists a dense open subset $U$ of the complement of the union of lines on $S$, consisting of smooth points distinct from $a$, such that, for every $q\in U$, the unique hyperplane containing $\mathbb T_qS$ and $a$ cuts an irreducible anticanonical curve $C(q,a)$, with an ordinary node at $q$, smooth elsewhere, avoiding $K\cup\Sing S$, and with tangent at $a$ different from $d$.
\end{lemma}
\begin{proof}
We blow up the minimal resolution of $S$ at $a$ and write $e_a$ for the new exceptional curve and $H$ for the pullback of the hyperplane class. Since the resolution of the nodes is crepant, the resulting surface has anticanonical divisor $H-e_a$, with
\[
 (H-e_a)^2=H^2+e_a^2=4-1=3.
\]
This divisor is nef: its intersection with the strict transform of a curve $D\subset S$ is $\deg D-\operatorname{mult}_aD\ge0$, with each exceptional curve over a node it is zero, and with $e_a$ it is one. The anticanonical linear system defines a birational morphism onto a normal cubic surface $\overline S\subset\PP^3$ with Du Val singularities, see \cite[Section~2.3, Proposition~2.6]{BH}. Since $(H-e_a)\cdot e_a=1$, this morphism restricts to an isomorphism from $e_a$ onto a line $\ell_a\subset\overline S$.

The rational map $S\dashrightarrow\overline S$ is the projection from $a$, so hyperplanes in $\PP^4$ through $a$ correspond to planes in $\PP^3$. Wherever this projection is a local isomorphism, a hyperplane containing $\mathbb T_qS$ and $a$ corresponds to the tangent plane to $\overline S$ at the image of $q$. We can therefore construct the required curve $C(q,a)$ using tangent-plane sections of $\overline S$.

Let $\overline U\subset\overline S_{\mathrm{reg}}$ be the open subset on which the second fundamental form is nondegenerate. We first show that $\overline U$ is nonempty. In characteristic zero, the closures of general fibers of the Gauss map are linear spaces \cite[Sections~1--2]{MT}. A normal cubic surface with Du Val singularities has finitely many lines, see \cite[Section~9.2.2, Table~9.1]{Dolgachev}, and is not a plane, so its Gauss map is generically finite. It is therefore generically unramified. Its differential is the second fundamental form, which consequently has rank two on a nonempty open subset. Thus $\overline U\ne\varnothing$.

We choose $\bar q$ in $\overline U$ away from the finitely many lines on $\overline S$. The tangent-plane section $\overline S\cap\mathbb T_{\bar q}\overline S$ then has an ordinary node at $\bar q$. If the plane cubic were reducible, a line component would pass through $\bar q$; indeed, a singular point on a reducible plane cubic curve off its line components would have to be a singular point of an irreducible conic, which is impossible. Our choice thus makes the section irreducible. Its arithmetic genus is one, and the node contributes delta-invariant one, so it is a rational nodal plane cubic. The same tangent-section construction is used in \cite[Section~5]{HT}.

No fixed point of $\PP^3$ belongs to every tangent plane of $\overline S$. Indeed, otherwise the polar derivative of its cubic equation in that point's direction would vanish on $\overline S$. This polynomial has degree two, so it would vanish identically, making $\overline S$ a cone in characteristic zero. A normal cubic cone has a vertex that is not Du Val. Thus, for any fixed point $r\in\PP^3$, the condition $r\notin\mathbb T_{\bar q}\overline S$ holds on a nonempty open subset of $\overline U$.

We now shrink $\overline U$ so that, for every $\bar q\in\overline U$, the projective tangent plane $\mathbb T_{\bar q}\overline S$ avoids the images of $K$, the singular points of $\overline S$, and the point of $\ell_a$ corresponding to $d$. Since $\overline U$ is irreducible, these finitely many open conditions leave a nonempty open subset. The tangent plane $\mathbb T_{\bar q}\overline S$ does not contain $\ell_a$, since the irreducible section $\mathbb T_{\bar q}\overline S\cap\bar S$ contains no line. This section meets $\ell_a$ transversely in one point; indeed, the plane cuts the line in degree one, at a smooth point of the section.

We let $Y$ denote the surface obtained above by resolving $S$ and blowing up $a$, with morphisms $\pi:Y\to S$ and $\phi:Y\to\overline S$. Since the plane section $\overline C=\mathbb T_{\bar q}\overline S\cap\overline S$ avoids the singularities of $\overline S$, its strict transform $\widetilde C\subset Y$ maps isomorphically to $\overline C$ under $\phi$. As $\phi$ maps $e_a$ isomorphically onto $\ell_a$, the curve $\widetilde C$ meets $e_a$ transversely at a single point $t$. Under the identification $e_a\simeq\PP(T_aS)$, this point represents a direction different from $d$. The image after blowing down $e_a$ is therefore a curve smooth at $a$, with tangent direction represented by $t$. The exceptional curves over the nodes of $S$ are disjoint from $\widetilde C$, so their contraction does not affect the curve. Thus $\pi(\widetilde C)$ is the required hyperplane section of $S$, smooth at $a$ with tangent different from $d$.

 Shrinking $\overline U$ further, we may assume that projection from $a$ induces an isomorphism from an open subset of $S\setminus\{a\}$ onto $\overline U$. Let $q$ be the preimage of $\bar q$; after a further shrinking, we may assume that $q$ lies on no line in $S$. Then $a\notin\mathbb T_qS$. Indeed, otherwise the line $aq$ would be tangent at $q$ to each of the two quadrics defining $S$ and would also meet each quadric at $a\ne q$; hence it would be contained in both quadrics, contrary to the choice of $q$. Consequently, the plane $\mathbb T_qS$ and the point $a$ span a projective $3$-space in $\PP^4$, which is the unique hyperplane containing them. These conditions define a nonempty open subset $U$ of $S$. Since $S$ is integral, $U$ is dense, proving the lemma.
\end{proof}

The following lemma extends the auxiliary point to a section of $f:\cS\to B$ so that the conclusions of Lemma~\ref{lem:movingcenter} hold at each $s_j(u)$ for all $u$ near $b$, together with the required condition at $z$.
\begin{lemma}\label{lem:auxiliarysection}
After shrinking $B$ around $b$, there exists a section $q:B\to\cS$, disjoint from the marked sections, such that for every $u\in B$ the point $q(u)$ lies on no line contained in $\cS_u$, and the following properties hold.
\begin{enumerate}
\item For every $u\in B$ and each $j=1,\ldots,16$, the tangent section $C(q(u),s_j(u))$ from Lemma \ref{lem:movingcenter} is irreducible, has an ordinary node at $q(u)$, is smooth elsewhere, and avoids the singularities of $\cS_u$ and the points $s_k(u)$ for $k\ne j$.
\item If $z$ is not exceptional for $\beta$, then $q(b)\ne a=\beta(z)$, and the curve $C(q(b),a)$ is irreducible, has an ordinary node at $q(b)$, is smooth elsewhere, and avoids the singularities and all sixteen points $s_k(b)$.
\item If $a=\beta(z)=s_i(b)$, the tangent direction of $C(q(b),s_i(b))$ at $s_i(b)$ differs from the direction represented by $z$.
\end{enumerate}
\end{lemma}
\begin{proof}
We apply Lemma~\ref{lem:movingcenter} at each $s_j(b)$, excluding the fifteen points $s_k(b)$ with $k\ne j$. If $z$ is not exceptional, we apply it also at $a=\beta(z)$, excluding all sixteen $s_k(b)$. If $z$ lies over $s_i(b)$, it specifies the direction at that point. Where $z$ prescribes no direction, we may use any $d_x\in\PP(T_x\cS_b)$ in Lemma~\ref{lem:movingcenter}. For each application the permissible auxiliary points form a nonempty open subset of $\cS_b$. Since $\cS_b$ is integral, we can choose $q_b$ in their intersection.

To extend $q_b$ to a section, we work on the relative smooth locus $\cS^{\sm/B}=\{x\in\cS:f\text{ is smooth at }x\}$. Each fiber of this smooth morphism is the smooth locus of a smooth or Du Val del Pezzo surface, and is strongly rationally connected by Xu's theorem \cite[Theorem~1.4]{Xu}. Recall that a smooth variety $V$ is \emph{strongly rationally connected} if every point of $V$ lies on a \emph{very free rational curve}, that is, the image of a morphism $h:\PP^1\to V$ for which $h^*T_V$ is ample. Equivalently, in the splitting $h^*T_V\simeq\bigoplus_i\OO_{\PP^1}(a_i)$, all $a_i$ are positive. Here $V=(\cS_u)_{\mathrm{reg}}$, so these curves lie entirely in the smooth locus of the fiber.

We use the following consequence of Hassett--Tschinkel's approximation theorem \cite[Theorem~15]{HT}: a smooth morphism over a smooth curve, with strongly rationally connected fibers and an existing section, admits a section through any prescribed point of a fiber. Here any one of the $s_j$ supplies the required section of $\cS^{\sm/B}\to B$. The theorem therefore ensures the existence of a section $q:B\to\cS^{\sm/B}$ with $q(b)=q_b$.

We now shrink $B$ around $b$ to preserve the required avoidance conditions. The union of the marked sections is closed in $\cS$. The lines contained in the fibers are parametrized by a projective relative Hilbert scheme; its universal family maps properly to $\cS$, with image equal to the union of these lines. This union is therefore closed as well. Since $q(b)=q_b$ belongs to neither closed set, their inverse images under $q$ can be removed from $B$.

For each $j$, the requirements that $C(q(u),s_j(u))$ is irreducible, has an ordinary node at $q(u)$ and no other singularities, and avoids the singularities of $\cS_u$ and the points $s_k(u)$ for $k\ne j$ are open conditions on $u$. They hold at $u=b$ by the choice of $q_b$. After shrinking $B$ once more, (1) holds for all $j$ and all $u\in B$. Properties (2) and (3) concern only the fiber over $b$ and already hold by the choice of $q_b$.
\end{proof}

\subsection{The conic bundle and the exceptional curve}\label{sec:conicmodel}
\subsubsection*{Construction of the conic bundle}
\paragraph{The blow-up and its cubic model.}
We choose the auxiliary section $q$ as in Lemma~\ref{lem:auxiliarysection} and blow up its image:
\[
 \nu:\widetilde{\cS}=\Bl_{q(B)}\cS\longrightarrow\cS.
\]
The threefold $\widetilde{\cS}$ is smooth, since both $\cS$ and the center $q(B)$ are smooth. For $u\in B$, its fiber $\widetilde{\cS}_u$ is the blow-up of $\cS_u$ at $q(u)$. We first identify $\widetilde{\cS}_u$ with a cubic surface and then construct a pencil of conics on it.

After shrinking $B$, we can choose projective coordinates depending algebraically on $u$ in which $q(u)=[1:0:0:0:0]$. With $x'=(x_1,\ldots,x_4)$, the two quadrics defining $\cS_u$ have the form
\[
 Q^{(i)}=x_0L_i(x')+F_i(x'),\qquad i=0,1,
\]
where $L_i$ and $F_i$ are respectively linear and quadratic forms in the fiber coordinates $x'$, with coefficients varying algebraically with $u\in B$. Since $q(u)$ is smooth on $\cS_u$, the forms $L_0,L_1$ are independent. Projection from $q(u)$ sends $[x_0:x']$ to $[x']\in\PP^3$. After eliminating $x_0$, the image is the cubic surface
\[
 \mathcal V_u=\{L_0F_1-L_1F_0=0\}\subset\PP^3.
\]
These cubic surfaces form a family in $B\times\PP^3$.

The projection with center $q(u)$ extends to $\widetilde{\cS}_u$ and maps its exceptional curve $e_{q(u)}$ onto the line
\[
 \ell_u=\{L_0=L_1=0\}\subset\mathcal V_u.
\]
To construct the inverse, we parametrize the line through $q(u)$ in direction $[y]\in\PP^3$ by $[x_0:t y]$. After removing the factor $t$ from the two quadric equations, the equations on the blow-up are
\[
 x_0L_i(y)+tF_i(y)=0,\qquad i=0,1.
\]
For $[y]\in\mathcal V_u$, these equations determine a unique $[x_0:t]$ unless all four coefficients $L_0(y)$, $L_1(y)$, $F_0(y)$, and $F_1(y)$ vanish. Such simultaneous vanishing would mean that the entire line through $q(u)$ in direction $[y]$ lies in $\cS_u$, contrary to Lemma~\ref{lem:auxiliarysection}. The kernel coordinates define a regular inverse, including along $\ell_u$, where $t=0$. These constructions vary algebraically with $u$, identifying $\widetilde{\cS}$ with the family of cubic surfaces. Moreover, the cubic surface $\widetilde{\cS}_u$ is smooth along $\ell_u$: there its gradient is $F_1\,dL_0-F_0\,dL_1$, which is nonzero by the same exclusion. On $\widetilde{\cS}_u$, the projection is given by the anticanonical system $|H-e_{q(u)}|$, where $H$ is the pullback of the hyperplane class of $\cS_u$.

\paragraph{The pencil of residual conics.}
A plane containing $\ell_u$ cuts the cubic surface in $\ell_u$ and a residual plane conic. For example, in the plane $L_1=\lambda L_0$, the cubic equation becomes
\[
 L_0(F_1-\lambda F_0)=0,
\]
so the residual conic is given by $F_1-\lambda F_0=0$ in that plane. The pencil of planes containing $\ell_u$ is parametrized by a projective line $T_u$. As $u$ varies, these lines form a $\PP^1$-bundle $T\to B$.

The pencil of residual conics has no base point. Indeed, outside $\ell_u$, a point belongs to only one plane of the pencil and hence cannot belong to every residual conic. A base point on $\ell_u$ would satisfy $F_1-\lambda F_0=0$ for every $\lambda$, and thus $F_0=F_1=0$, which we have excluded. This base-point-free pencil therefore defines a morphism to $T_u$, and these fiberwise morphisms fit together into
\[
 h_q:\widetilde{\cS}\longrightarrow T,\qquad T\longrightarrow B.
\]
The fiber of $h_q$ over a point of $T_u$ is the corresponding residual conic on $\widetilde{\cS}_u$. Each plane of the pencil cuts out $\ell_u$ and a unique residual conic. Every such conic is geometrically connected: it is a smooth conic, two intersecting lines, or a double line. Thus a finite cover of degree greater than one in the Stein factorization would disconnect a general fiber, which is impossible. Since $T$ is normal, the remaining finite birational morphism is an isomorphism. Since the cubic is integral, no plane of the pencil is a component of it, and every residual conic has pure dimension one. Both $\widetilde{\cS}$ and $T$ are smooth, so $h_q$ is flat by miracle flatness \cite[Lemma~10.128.1, Tag~00R4]{Stacks}. The composition $\widetilde{\cS}\to T\to B$ is the original surface fibration after blowing up $q(B)$.

\paragraph{The correspondence with nodal tangent sections.}
Recall that $q$ is the auxiliary section constructed in Lemma~\ref{lem:auxiliarysection}. For each curve $C(q(u),a)$ in that lemma, its strict transform $D_{u,a}$ under the blow-up $\nu_u:\widetilde{\cS}_u\to\cS_u$ is the corresponding residual conic. Conversely, $C(q(u),a)=\nu_u(D_{u,a})$ is the proper transform of this conic on $\cS_u$ under the inverse birational correspondence.

Indeed, the hyperplane defining $C(q(u),a)$ projects to a plane containing $\ell_u$, and its section of the cubic is $\ell_u+D_{u,a}$. The class of the residual conic is therefore
\[
 [D_{u,a}]=(H-e_{q(u)})-e_{q(u)}=H-2e_{q(u)}.
\]
This is also the class of the strict transform of $C(q(u),a)$, whose multiplicity at $q(u)$ is two. The two branches meet the exceptional curve over $q(u)$ at distinct points. Since $C(q(u),a)$ is smooth elsewhere, $D_{u,a}$ is a smooth conic and $\nu_u|_{D_{u,a}}:D_{u,a}\to C(q(u),a)$ is its normalization map.

The conic through the lift of $s_j(u)$ is smooth, since it is the normalization of $C(q(u),s_j(u))$. For distinct $j,k$, these conics have distinct parameters in $T_u$, since $C(q(u),s_j(u))$ avoids $s_k(u)$ by Lemma~\ref{lem:auxiliarysection}(1). If $z$ is nonexceptional for $\beta$, the conic through the lift of $a=\beta(z)$ is also smooth. Its parameter in $T_b$ differs from those of the sixteen conics through the lifts of the $s_j(b)$, since $C(q(b),a)$ avoids all sixteen points by part~(2) of the same lemma. Finally, the parameters of singular conics form a closed subset of $T_b\simeq\PP^1$. It omits the parameters of the smooth conics just described and is therefore finite.

\begin{lemma}\label{lem:conicbundle}
For the auxiliary section $q$ of Lemma~\ref{lem:auxiliarysection}, the $\PP^1$-bundle $T\to B$ and morphism
\[
 h_q:\widetilde{\cS}=\Bl_{q(B)}\cS\longrightarrow T
\]
are as constructed above. The morphism is flat with plane-conic fibers and has the following properties.
\begin{enumerate}
\item For every $u\in B$, the conics through the lifts of the points $s_j(u)$ are smooth and have pairwise distinct parameters in $T_u$. They are the strict transforms of $C(q(u),s_j(u))$, and the blow-down restricts to their normalization maps.
\item If $z$ is nonexceptional for $\beta$, the conic through the lift of $a=\beta(z)$ is smooth and has a parameter in $T_b$ different from the parameters of the sixteen conics in (1) over $b$. It is the strict transform of $C(q(b),a)$.
\item Only finitely many parameters in $T_b$ correspond to singular conics.
\end{enumerate}
\end{lemma}
\begin{proof}
The construction above proves flatness and identifies the fibers as residual plane conics. Their stated properties follow from the correspondence with the tangent sections in Lemma~\ref{lem:auxiliarysection}, as established in the preceding paragraph.
\end{proof}

\subsubsection*{Retaining the exceptional curve containing $z$}
When $z$ is nonexceptional we use $\mathscr M=\widetilde{\cS}$. If $z$ lies above $s_i(b)$, we retain that blow-up and set $\mathscr M=\Bl_{s_i(B)}\widetilde{\cS}$. The image of $s_i(B)$ in $T$ is a section of $T\to B$, hence a Cartier divisor, and is disjoint from the discriminant of $h_q$. For each $u\in B$, the smooth conic through the lift of $s_i(u)$ becomes
\begin{equation}\label{eq:splitconic}
 e_i\cup C',
\end{equation}
and the two reduced projective lines meet transversely. Here $e_i$ is the exceptional curve and $C'$ the strict transform of the smooth conic. On the fiber over $b$, let $t_i=e_i\cap C'$ be the unique intersection point of the two components. Under the identification $e_i\simeq\PP(T_{s_i(b)}\cS_b)$, the point $t_i$ represents the tangent direction of $C(q(b),s_i(b))$ at $s_i(b)$. Indeed, the strict transform of a smooth curve through $s_i(b)$ meets $e_i=\PP(T_{s_i(b)}\cS_b)$ at its tangent direction. Property (3) of Lemma~\ref{lem:auxiliarysection} requires this direction to differ from the direction represented by $z\in e_i$. Hence $t_i\ne z$, so $z\in e_i\setminus C'$ is a smooth point of the reducible conic.

\subsubsection*{Embedding the conic family in a projective-plane bundle}
We now place all fibers of $h:\mathscr M\to T$ in projective planes varying algebraically with the parameter in $T$. This includes the reducible fibers introduced by the additional blow-up. Our goal is to describe the whole family by a quadratic equation in three homogeneous coordinates, locally on $T$.

\paragraph{The line bundle on the fibers.}
Both $\mathscr M$ and $T$ are smooth. Since every fiber of $h:\mathscr M\to T$ has pure dimension one, $h$ is flat by miracle flatness. Each fiber is a Gorenstein curve, since it is either a plane conic or two smooth rational curves meeting transversely at one point. Consequently, the relative dualizing sheaf $\omega_{\mathscr M/T}$ is invertible and restricts on each fiber $C$ to its dualizing sheaf $\omega_C$, see \cite[Lemma~109.19.3, Tag~0E6R]{Stacks}. We use its inverse
\[
 \mathcal L_h=\omega_{\mathscr M/T}^{-1}.
\]
We shall show that $\mathcal L_h|_C$ has three independent sections and embeds $C$ as a plane conic.

For a fiber $C$ unchanged by the additional blow-up, its plane embedding is already known. For a degree-two curve in $\PP^2$, adjunction implies $\omega_C\simeq\OO_C(-1)$, hence $\mathcal L_h|_C\simeq\OO_C(1)$. This remains valid for a reducible conic or a double line. If $F=0$ is its quadratic equation, the exact sequence
\[
 0\longrightarrow\OO_{\PP^2}(-1)
 \xrightarrow{\ \cdot F\ }\OO_{\PP^2}(1)
 \longrightarrow\OO_C(1)\longrightarrow0
\]
shows that the three linear coordinate functions restrict to a basis of $H^0(C,\mathcal L_h|_C)$ and that $H^1(C,\mathcal L_h|_C)=0$. These sections have no common zero and define the original embedding of $C$ into $\PP^2$.

On a new fiber $C=e_i\cup C'$ from \eqref{eq:splitconic}, the two components meet at $t=e_i\cap C'$. On each component the dualizing sheaf of $C$ is the canonical bundle of that component with the intersection point added. Thus
\[
 \omega_C|_{e_i}\simeq\omega_{e_i}(t)\simeq\OO_{\PP^1}(-1),
 \qquad
 \omega_C|_{C'}\simeq\omega_{C'}(t)\simeq\OO_{\PP^1}(-1).
\]
It follows that $\mathcal L_h$ has degree one on each component. A section on $C$ is a pair of sections on the components whose values agree at $t$. Equivalently, there is an exact sequence
\[
 0\longrightarrow\mathcal L_h|_C
 \longrightarrow\mathcal L_h|_{e_i}\oplus\mathcal L_h|_{C'}
 \longrightarrow\mathcal L_h|_t\longrightarrow0,
\]
where the last map takes the difference of the two values at $t$. Each component contributes two independent sections, and matching their values imposes one condition. Hence again $h^0(C,\mathcal L_h|_C)=3$ and $h^1(C,\mathcal L_h|_C)=0$. We can take as a basis one section nonzero at $t$ and two others vanishing there, one supported on each component. The resulting map embeds the two components as two distinct lines in $\PP^2$, meeting at the image of $t$.

\paragraph{Putting the fiberwise embeddings together.}
We write $E_h=h_*\mathcal L_h$ for the direct image. The computations above show that the space of sections on every fiber has dimension three and that its first cohomology vanishes. The cohomology and base change lemma \cite[Lemma~30.22.1, Tag~07VK]{Stacks} therefore implies that $E_h$ is locally free of rank three and that, for every $\tau\in T$,
\[
 E_h\otimes k(\tau)\simeq H^0(h^{-1}(\tau),\mathcal L_h|_{h^{-1}(\tau)}).
\]
In particular, a local basis of $E_h$ corresponds to three sections of $\mathcal L_h$ whose restrictions form a basis on each fiber. Since these sections have no common zero on any fiber, the evaluation homomorphism is surjective:
\[
 h^*E_h\longrightarrow\mathcal L_h.
\]
Explicitly, this map evaluates a section at a point of $\mathscr M$. In a local trivialization of $E_h$, with basis $s_0,s_1,s_2$, it defines the map $x\mapsto[s_0(x):s_1(x):s_2(x)]$. Changes of basis change these coordinates by a projective linear transformation. Thus the local maps fit together into a morphism to the projective-plane bundle associated with $V=E_h^\vee$, where $\PP_T(V)$ parametrizes one-dimensional subspaces of the fibers of $V$. Since $h$ is proper and these maps are closed embeddings on every fiber, the resulting morphism is a relative closed embedding:
\begin{equation}\label{eq:conicembedding}
 \mathscr M\hookrightarrow\PP_T(V),\qquad V=(h_*\omega_{\mathscr M/T}^{-1})^\vee.
\end{equation}

\paragraph{The quadratic equation.}
Locally on $T$, we choose the three homogeneous coordinates $w_0,w_1,w_2$ just described. The embedded family is cut out by one equation
\[
 Q(w_0,w_1,w_2)=\sum_{0\le i\le j\le2}a_{ij}w_iw_j=0,
\]
with coefficients $a_{ij}$ regular on that open subset of $T$. To see how these equations fit together, we use the multiplication map
\[
 \operatorname{Sym}^2E_h\longrightarrow h_*(\mathcal L_h^{\otimes2}).
\]
On each plane conic, the six quadratic monomials span the five-dimensional space $H^0(C,\OO_C(2))$; their kernel is the one-dimensional space generated by its defining equation. The cohomology and base change lemma \cite[Lemma~30.22.1, Tag~07VK]{Stacks} implies that $h_*(\mathcal L_h^{\otimes2})$ is locally free of rank five. The multiplication map is surjective on every fiber, so its kernel is a line bundle $J$ on $T$. Its inclusion in $\operatorname{Sym}^2E_h$ is equivalently a quadratic form
\[
 Q\in H^0\bigl(T,\operatorname{Sym}^2V^\vee\otimes L\bigr),
 \qquad L=J^\vee.
\]
For $w\in V_\tau$, the value $Q(w)$ lies in $L_\tau$. On overlaps, its local equations differ by an invertible factor, together with the change of homogeneous coordinates, so their zero sets define the same family $\mathscr M$.

\begin{lemma}\label{lem:conicbundle-embedding}
The morphism $h:\mathscr M\to T$ is a conic bundle: it is flat and projective, and every geometric fiber is a plane conic. The sheaf
\[
 V=(h_*\omega_{\mathscr M/T}^{-1})^\vee
\]
is locally free of rank three. There exist a line bundle $L$ on $T$ and a quadratic form
\[
 Q\in H^0(T,\operatorname{Sym}^2V^\vee\otimes L)
\]
such that the embedding \eqref{eq:conicembedding} identifies $\mathscr M$ with the zero scheme $\{Q=0\}\subset\PP_T(V)$. This description includes the reducible fibers $e_i\cup C'$ introduced by the retained blow-up.
\end{lemma}
\begin{proof}
Flatness and the description of the fibers were established above. The relative anticanonical embedding \eqref{eq:conicembedding} proves projectivity. If $J$ is the line bundle of quadratic relations constructed above, then $L=J^\vee$ and $Q$ is the associated quadratic form.
\end{proof}

\subsection{Reflections over the parameter curve}\label{sec:reflectioninterpolation}
We use the quadratic equation $Q=0$ of Lemma~\ref{lem:conicbundle-embedding} to construct a birational involution of $\mathscr M$ that preserves each fiber of $h:\mathscr M\to T$. On the conic containing the point represented by $z$, it must be regular and move that point to one whose image in $\cS_b$ lies on no line and differs from every $s_j(b)$. If $z$ lies on the retained exceptional curve, it must in particular move $z$ to the other component of the reducible conic.

We first describe reflections on individual conics. For each section $s_j$ still to be blown up, let $c_j:B\to\mathscr M$ be its lift. We prescribe reflections on the conic through $z$ and on the singular conics over $b$. On the conic through each $c_j(b)$, we choose a reflection that moves this point. Interpolating these choices will allow us to control indeterminacy on $Z_B$ in Subsection~\ref{sec:blowupcomparison}.

\paragraph{A reflection on one conic.}
On a smooth conic $C=\{Q=0\}\subset\PP^2$, we choose the center $p$ in $\PP^2\setminus C$. Each line through $p$ meets $C$ in two points counted with multiplicity. Interchanging these points defines an involution of $C$, called a reflection with center $p$. The tangency points are fixed. Equivalently, projection from $p$ is a degree-two morphism from $C$ to the projective line parametrizing lines through $p$, and the reflection exchanges its two sheets.

To express this construction in coordinates, let $B_Q$ be the symmetric bilinear form associated with $Q$, normalized by $B_Q(w,w)=Q(w)$. A representative $v$ of $p$ satisfies $Q(v)\ne0$. For $[w]\in C$, the equation on the line joining $[w]$ to $p$ is
\[
 Q(w+tv)=2tB_Q(w,v)+t^2Q(v).
\]
Besides the intersection $t=0$, its second intersection with $C$ is therefore represented by the vector
\begin{equation}\label{eq:reflection}
 R_v(w)=Q(v)w-2B_Q(w,v)v.
\end{equation}
The identities
\begin{equation}\label{eq:reflectionidentities}
 Q(R_v(w))=Q(v)^2Q(w),\qquad R_v^2=Q(v)^2I,\qquad
 \det R_v=-Q(v)^3
\end{equation}
are polynomial identities and do not require nonsingularity of $Q$. Where $Q(v)$ is a unit, the induced projective map is an automorphism of the whole conic family.

A vector $v$ with $Q(v)\ne0$ will be called \emph{anisotropic}; geometrically, this means that its projective point $[v]$ lies outside the conic. For such a vector the reflection is invertible, by \eqref{eq:reflectionidentities}. On a smooth conic it can move a given point to any distinct point: for isotropic representatives $w,w'$ of the two points, we have $B_Q(w,w')\ne0$, and the reflection with center $[w-w']$ interchanges them. For a conic consisting of two distinct lines, the same geometric construction applies when the center lies outside both components. A line through the center meets the two components in points which are interchanged; the common node is fixed. The invertible linear map \eqref{eq:reflection} defines the involution at that node as well. On the split conic $XY=0$, the formula for a vector $(A,B,D)$ is
\[
 [X:Y:Z]\longmapsto[-A^2Y:-B^2X:ABZ-D(BX+AY)].
\]
If $AB\ne0$, it exchanges the components; $(A,B,D)=(1,1,k+\ell)$ takes $[0:1:k]$ to $[1:0:\ell]$. In the split case we can therefore send $z$ to any point of the other component except the node. In the smooth case we can send $z$ to any other point of its conic. The component on which we choose the reflected point is the strict transform of the irreducible tangent section $C(q,\beta(z))$ and maps onto it by normalization. The tangent quartic meets the union of the lines contained in the del Pezzo fiber $\cS_b$, together with the marked points, in only finitely many points. The strict transform also meets the auxiliary exceptional divisor in finitely many points. We choose the reflected point away from these intersections and from the preimages of the excluded points on $\cS_b$.

\paragraph{Varying the center with the conic.}
For the family $\mathscr M\subset\PP_T(V)$, the centers must vary algebraically with the conic parameter. Here $V$ has rank three, so each fiber $\PP(V_\tau)$ is a projective plane. A section $p$ of this projective-plane bundle chooses a center $p(\tau)$ in each plane. Wherever $p(\tau)$ lies outside the conic $h^{-1}(\tau)$, the reflection is an automorphism of that conic.

A local representative of $p$ is a nonzero vector $v$. Replacing $v$ by $\lambda v$, with $\lambda$ invertible, multiplies $R_v$ by $\lambda^2$, so its projective action does not change. We may therefore represent the centers by a section of $V\otimes A$, where $A$ is a line bundle on $T$, provided that the section is nonzero on the locus under consideration. This freedom to twist will allow us to prescribe the centers at finitely many parameters. We shall require those centers to lie outside their conics; at other parameters, the center may lie on the conic, and the resulting map is then only rational. Its possible indeterminacy will be examined below.

\paragraph{Choosing the reflection on the other conics.}
We must now choose the vector $v$ algebraically on the complete curve $T_b\simeq\PP^1$, which parametrizes all conics in the surface over $b$. Let $t_0\in T_b$ be the parameter of the conic through $z$, and let $\Gamma$ be the finite set of parameters of singular conics. For the lifts $c_j$ defined above, Lemma~\ref{lem:conicbundle} shows that the parameters $a_j=h(c_j(b))$ are pairwise distinct, lie outside $\Gamma$, and differ from $t_0$. At $t_0$ we prescribe a vector defining the reflection just chosen. At every point of $\Gamma\cup\{a_j\}$ we prescribe an anisotropic vector; at $a_j$ we require also that its reflection move $c_j(b)$. Since the reflection preserves $h$ and the $a_j$ are distinct, the image of $c_j(b)$ also differs from every $c_k(b)$ with $k\ne j$.

The condition $Q(v)\ne0$ makes the reflection an automorphism on the whole conic at that parameter. In particular, a preimage of $c_j(b)$ is then a single point, and no such conic is contracted to $c_j(b)$.

\paragraph{Interpolating the prescribed values.}
Let $D_*$ be the reduced divisor supported on $\{t_0\}\cup\Gamma\cup\{a_j\}$. We choose a twist $A$ sufficiently positive on the fibers of $T\to B$ so that the bundle
\[
 F_b=V|_{T_b}\otimes A|_{T_b}
\]
satisfies $H^1(T_b,F_b(-D_*))=0$ and $F_b(-D_*)$ is globally generated. We keep this twist fixed as the section varies. The vanishing implies that evaluation
\[
 H^0(T_b,F_b)\longrightarrow H^0(D_*,F_b|_{D_*})
\]
is surjective. Thus the sections taking the prescribed values form a nonempty affine space $W$, a translate of $H^0(T_b,F_b(-D_*))$.

For each $t\in T_b\setminus D_*$, global generation of $F_b(-D_*)$ implies that the condition $v(t)=0$ has codimension three in $W$. The incidence of such pairs $(t,v)$ in $(T_b\setminus D_*)\times W$ therefore has dimension at most $\dim W-2$, as does the closure of its image in $W$. A general section $v_0\in W$ has no zero outside $D_*$. Its prescribed values on $D_*$ are nonzero, so it is nowhere zero on $T_b$.

The same vanishing implies $H^1(T_b,F_b)=0$. By cohomology and base change, after shrinking $B$ around $b$, the sections of $V\otimes A$ on the fibers form a vector bundle whose fiber at $b$ is $H^0(T_b,F_b)$. The chosen vector $v_0$ extends to a local section of this bundle, hence to a section $v$ of $V\otimes A$ over $T$.

If $Q$ takes values in a line bundle $L$ and $v$ is a section of $V\otimes A$, then
\[
 R_v:V\longrightarrow V\otimes L\otimes A^2.
\]
Its projectivization is a rational self-map of $\PP_T(V)$; tensoring a vector bundle by a line bundle does not change its projective bundle.

\paragraph{The points where the reflection may be undefined.}
Let $\rho$ be the resulting rational involution of $\mathscr M$ and let $n=Q(v)$. On $T_b$, $n$ is not identically zero. Its zeros occur only at smooth conics, where $v$ is nonzero isotropic. At such a zero, \eqref{eq:reflection} is undefined on the conic only where $B_Q(w,v)=0$, the tangent line at $[v]$. It has just one point of support on the smooth conic. Every singular conic lies in the automorphism locus. Therefore
\begin{equation}\label{eq:rhofinite}
 \Indet(\rho)\cap\mathscr M_b\text{ is finite}.
\end{equation}
The involution is an isomorphism near each $c_j(b)$ and sends it outside the set of points $c_j(b)$. After shrinking $B$, $\rho$ is an isomorphism on a neighborhood of each $c_j(B)$, and $\rho(c_j(B))$ is disjoint from every $c_k(B)$.

\begin{lemma}\label{lem:reflection-on-conic-model}
Let $z_{\mathscr M}$ be the point of $\mathscr M_b$ corresponding to $z$, and let $\mu:\mathscr M\to\cS$ be the natural morphism. After shrinking $B$ around $b$, there exists a birational involution $\rho$ of $\mathscr M$ with $h\rho=h$ and the following properties:
\begin{enumerate}
\item The set $\Indet(\rho)\cap\mathscr M_b$ is finite.
\item The map $\rho$ is an isomorphism between neighborhoods of $z_{\mathscr M}$ and $\rho(z_{\mathscr M})$. Its image $\rho(z_{\mathscr M})$ avoids the auxiliary exceptional divisor, and $\mu(\rho(z_{\mathscr M}))$ is a smooth point of $\cS_b$ lying on no line and distinct from all marked points.
\item For each of the sections $c_j:B\to\mathscr M$ defined above, the map $\rho$ is an automorphism on a neighborhood of the conic through $c_j(b)$. It is also an isomorphism on a neighborhood of $c_j(B)$, and $\rho(c_j(B))$ is disjoint from every $c_k(B)$.
\end{enumerate}
\end{lemma}
\begin{proof}
The prescribed reflection at $t_0$ has the properties in (2). At each $a_j$, the prescribed vector is anisotropic and its reflection moves $c_j(b)$; since the parameters $a_j$ are distinct, (3) holds after a further shrinking of $B$. The interpolation constructs these reflections in a single rational involution, whose indeterminacy on $\mathscr M_b$ is finite by \eqref{eq:rhofinite}.
\end{proof}

\subsection{Comparing the blow-ups}\label{sec:blowupcomparison}
The involution in Lemma~\ref{lem:reflection-on-conic-model} acts on $\mathscr M$, whereas Theorem~\ref{thm:relocation} requires one on $Z_B$. Transferring a rational map between birational models can introduce additional indeterminacy. We must therefore check what happens when we blow up the curves $c_j(B)$ and contract the auxiliary exceptional divisor. The following comparison also verifies that the transferred map remains a local isomorphism at $z$ and at its image.

If $\beta(z)=s_i(b)$, the conic model already contains the exceptional divisor over $s_i(B)$. In either case the sections $c_j$ are exactly those still to be blown up. To compare it with $Z_B$, we blow up the auxiliary section in $Z_B$. The centers are disjoint by Lemma~\ref{lem:auxiliarysection}, so we may interchange the blow-up along $q$ with those along the $c_j$. The common model is therefore
\[
 \mathscr Y=\Bl_qZ_B=\Bl_{\bigcup c_j(B)}\mathscr M,
 \qquad \pi:\mathscr Y\to Z_B,\qquad\alpha:\mathscr Y\to \mathscr M.
\]
On this common model we transport $\rho$ to the rational self-map $\theta:Z_B\dashrightarrow Z_B$ defined by
\begin{equation}\label{eq:theta}
 \theta=\pi\alpha^{-1}\rho\alpha\pi^{-1}.
\end{equation}
The map $\pi$ contracts only the auxiliary exceptional divisor, while $\alpha$ contracts the exceptional divisors over the curves $c_j(B)$. Formula~\eqref{eq:theta} conjugates $\rho$ by the birational comparison $\alpha\pi^{-1}$, so it is an involution preserving $B$. We now check indeterminacy on $(Z_B)_b$ and the local isomorphism property at $z$.

\paragraph{Indeterminacy on the fiber over $b$.}

Away from the single point $q(b)$, $\pi^{-1}$ is a local isomorphism. If a point lies on an $\alpha$-exceptional curve, its image is $c_j(b)$. The map $\rho$ is regular there and its image avoids the points $c_j(b)$. Thus $\alpha^{-1}$ is regular near that image and the final $\pi$ is a morphism. The composition is regular near \emph{every} point of this exceptional curve.

Away from the $\alpha$-exceptional curves, $\alpha$ is a local isomorphism. The only possible failures are a base point of $\rho$, or a point whose $\rho$-image is some $c_j(b)$. The first set is finite by \eqref{eq:rhofinite}. In the second case the conic parameter is $a_j$, where $n$ is a unit. The reflection is therefore an automorphism on a neighborhood of the complete conic, and the preimage of $c_j(b)$ is one point. Consequently
\begin{equation}\label{eq:bound}
 \#\bigl(\Indet(\theta)\cap(Z_B)_b\bigr)
 \le 1+\#\{a\in T_b:n(a)=0\}+\#\{j\}.
\end{equation}
Since $\theta$ is an involution, \eqref{eq:bound} also proves \eqref{eq:finiteindet}.

\paragraph{The image of $z$.}
Near $z$, the comparison morphisms are local isomorphisms: its exceptional divisor, when present, was retained in $\mathscr M$. The reflection is an automorphism near the conic through $z_{\mathscr M}$. By Lemma~\ref{lem:reflection-on-conic-model}(2), the point $\rho(z_{\mathscr M})$ avoids the curves $c_j(B)$ and the auxiliary exceptional divisor. Hence $\theta$ is a local isomorphism at $z$ and moves it to the required point $y$. The resulting involution on the original model has the following properties.

\begin{lemma}\label{lem:reflection-on-marked-model}
The map $\theta$ defined by \eqref{eq:theta} is a birational involution of $Z_B$ satisfying $g\theta=g$. It is an isomorphism between neighborhoods of $z$ and $y=\theta(z)$, where $\beta(y)$ is a smooth point of $\cS_b$ lying on no line and distinct from all marked points. Moreover, both
\[
 \Indet(\theta)\cap(Z_B)_b
 \quad\text{and}\quad
 \Indet(\theta^{-1})\cap(Z_B)_b
\]
are finite.
\end{lemma}
\begin{proof}
Birational conjugation preserves the involution property, and all maps in \eqref{eq:theta} preserve $B$. The assertions concerning $z$ and $y$ follow from the preceding local comparison and Lemma~\ref{lem:reflection-on-conic-model}. The bound \eqref{eq:bound} proves finiteness of the first indeterminacy set; the second is the same set because $\theta^{-1}=\theta$.
\end{proof}

\begin{proof}[Proof of Theorem~\ref{thm:relocation}]
For the point $z$ fixed at the beginning of this section, Lemma~\ref{lem:reflection-on-marked-model} constructs a birational involution $\theta$ of $Z_B$ preserving $g$. It is an isomorphism near $z$ and its image $y$, the point $\beta(y)$ lies on no line and is not marked, and the inverse indeterminacy on $(Z_B)_b$ is finite. These are the required properties.
\end{proof}

\section{Ellipticity of the small resolutions}\label{app:smallres}
For $C\in\mathcal U_0$, we complete the proof of assertion (3) of Theorem~\ref{thm:main} by establishing ellipticity of the two small resolutions.

\begin{theorem}\label{thm:smallres}
The small resolutions $X_R$ and $X_L$ are algebraically elliptic. They are nonrational.
\end{theorem}

A point of $X_R$ outside the sixteen exceptional curves has a unique lift to $Z$. Near that lift we identify the two base neighborhoods by $\beta_R$ and compose the local dominating sprays on $Z$ with $\beta_R$. The differential of $\beta_R$ at that lift is an isomorphism. The same argument applies to $X_L$.

At a point of an exceptional curve, this argument does not apply. We shall construct a birational involution taking that point to one where three independent tangent directions can be chosen along the two del Pezzo fibrations, and transfer sprays from its image. Here an additional difficulty arises: we need sprays along both fibrations, although only the right fibration is a morphism on $X_R$. Curves along the left fibration may pass outside the neighborhood of the right fiber used in Section~\ref{sec:relocation}. We therefore need a birational involution defined as a rational map on the whole threefold, with control of its indeterminacy on both surfaces through the image point. The complete curves used for the sprays must be disjoint from the indeterminacy locus of the inverse map.

We work with $X_R$; the same argument applies to $X_L$ after interchanging the two constructions. We first introduce the notation needed to state precisely what the involution must do.

Let
\[
 \beta_R:Z\longrightarrow X_R,\qquad
 C_R=\bigcup_{p\in\Sigma} C_p^R,\qquad C_p^R=\pi_R^{-1}(p).
\]
Thus $\beta_R$ is the blow-up of the sixteen disjoint sections $C_p^R$ of $\varphi_R$. On $X_R\setminus C_R$ it is an isomorphism. The left kernel defines a rational map
\[
 \psi_L:X_R\dashrightarrow \PP_L^1,\qquad \psi_L\circ\beta_R=g_L.
\]
Its indeterminacy locus is $C_R$. For $y\notin C_R$, we write
\[
 S_R(y)=\varphi_R^{-1}(\varphi_R(y)),
\]
and let $S_L(y)$ denote the closure in $X_R$ of the fiber of $\psi_L$ through $y$.

The next theorem is the main construction of this section. Its last assertion is needed when curves on $X_R$ are lifted to $Z$: a marked point in the indeterminacy locus would lift to an entire exceptional curve in the corresponding right fiber.

\begin{theorem}\label{thm:smallrelocation}
For every $x\in C_R$ there is a right-fiber-preserving birational involution
\[
 \theta:X_R\dashrightarrow X_R
\]
which is a local isomorphism at $x$. Its image $y=\theta(x)$ lies off $C_R$, and $\pi_R(y)\in X\setminus\Sigma$ lies on no projective line in either of the two del Pezzo surfaces in $\PP^4$ through it. Moreover, $\rk d(\varphi_R,\psi_L)_y=2$. If $J=\Indet(\theta^{-1})$, then
\[
 J\cap S_R(y)\qquad\text{and}\qquad J\cap S_L(y)
\]
are finite. Moreover $J\cap C_R$ is finite, and $J\cap S_R(y)$ misses the sixteen marked points.
\end{theorem}

We prove this theorem in three steps. We first show that a tangent section in a fixed right fiber meets varying left fibers. We then extend the reflection construction to a rational self-map of the whole threefold. Finally, we choose the reflection so that its indeterminacy has finite intersection with the left surface through its image point. The construction of sprays follows in the last subsection.

\subsection{Varying the left parameter on a tangent section}\label{app:oppositekernel}
On a fixed fiber $S_b=\varphi_R^{-1}(b)$ we write $H=\OO_{S_b}(1)$. The corresponding fiber of $g_R$ is the blow-up of $S_b$ at its sixteen marked points. We denote its exceptional curves by $E_1,\ldots,E_{16}$.

\begin{lemma}\label{lem:opposite-nonconstant}
The morphism to $\PP_L^1$ induced by $\psi_L$ on $\Bl_{16}S_b$ has fiber class
\[
 2H-\sum_{j=1}^{16}E_j.
\]
Let $D\subset S_b$ be an irreducible tangent hyperplane section as in Lemma~\ref{lem:movingcenter}, and suppose that $D$ contains at most one marked point. Then $\psi_L|_D$ is nonconstant.
\end{lemma}
\begin{proof}
On an affine left-kernel chart, the left kernel on $S_b$ is the ratio of two quadratic forms. Their common zero locus is the locus where both rows of $M$ vanish, hence where $M=0$. On the fixed right surface this is $\Sigma$. The four original quadrics meet transversely at each point of $\Sigma$ by Proposition~\ref{prop:geometry}, so these base points are simple. The pencil therefore has class $2H$ on $S_b$, and after blowing up the sixteen base points its fiber class is $2H-\sum_jE_j$.

The strict transform of $D$ has intersection number $8$ with the left fiber class if $D$ misses the marked points, and $8-\operatorname{mult}_{p_i}D\ge6$ if $p_i$ is its only marked point, so $\psi_L|_D$ is nonconstant.
\end{proof}

In particular, for $x\in C_p^R$, there is an irreducible tangent section $D_0(q,x)$ through $x$ avoiding the other fifteen marked points, see Lemma~\ref{lem:movingcenter}. The left parameter $\psi_L$ is nonconstant on this curve. Thus moving a point along the curve also varies the left surface through it.

\paragraph{Avoiding the lines of both fibrations.}
We may also choose $q$ so that $D_0(q,x)$ is not contained in the locus of lines in fibers of the left pencil. The relative Hilbert scheme of lines in the left pencil is proper with finite fibers, hence finite over $\PP_L^1$. Its universal family has dimension at most two, and its image in $X$ is closed. We identify $S_b$ with its image under $\pi_R$. If it were a surface component of the image of the universal family, a component of the universal family would cover $S_b$ by projective lines. This contradicts the finiteness of lines on $S_b$, both in the smooth and the two-nodal cases. Thus the intersection of the left line locus with $S_b$ is a proper closed subset of dimension at most one. A fixed irreducible section $D_0(q,x)$ determines $q$ as its unique node, so these sections cannot all be among the finitely many curve components of that intersection. Containment in the fixed closed locus is a closed condition in this flat family of hyperplane sections. A general $q$ therefore avoids containment, and a general point of $D_0(q,x)$ lies on no line of either del Pezzo fibration.

\paragraph{Avoiding the critical locus.}
The critical locus also meets $S_b$ in a finite set. On each coordinate plane, the two restricted diagonal equations have coefficient rank two: otherwise three columns would be proportional, contrary to condition~(ii). They determine at most one projective squared-coordinate point, with finitely many inverse images under coordinate squaring. Thus a general point of $D_0(q,x)$ lies in the rank-two locus.

\subsection{A global family of reflections on \texorpdfstring{$X_R$}{X\_R}}\label{app:globalreflections}
To accommodate the left-relative curves, we extend the reflection to the whole of $X_R$ and bound its indeterminacy on the right fibers outside the neighborhood used in Section~\ref{sec:relocation}.

For $x\in C_p^R$ and $b=\varphi_R(x)$, we choose the auxiliary section $q$ as in the proof of Theorem~\ref{thm:relocation}. We take $q(b)$ general as in Subsection~\ref{app:oppositekernel}. On $S_b$ the tangent hyperplane cuts out the curve $D_0(q(b),x)$ through $x$. For each $j$, the curve $C(q(b),s_j(b))$ has the properties in Lemma~\ref{lem:auxiliarysection}(1). Since $X_R$ is proper over $\PP_R^1$, $q$ extends over the whole base. We blow up its image:
\[
 \pi:\mathscr N=\Bl_{q(\PP^1)}X_R\longrightarrow X_R.
\]
Over a Zariski neighborhood $B^\circ\ni b$, the tangent-section pencil is a regular plane-conic bundle
\[
 h^\circ:\mathscr N^\circ\longrightarrow T^\circ,
\]
where $T^\circ\to B^\circ$ is a $\PP^1$-bundle. Let $V^\circ$ and $L^\circ$ be the rank-three bundle and line bundle for its relative anticanonical realization, with quadratic form
\[
 Q^\circ\in H^0\!\left(T^\circ,\operatorname{Sym}^2(V^\circ)^\vee\otimes L^\circ\right).
\]

\paragraph{Completing the parameter surface.}
We choose a smooth projective compactification $\overline T$ of $T^\circ$ for which the map to $\PP_R^1$ extends. A finite-presentation extension of $V^\circ$ to $\overline T$ exists by the extension theorem for coherent sheaves \cite[Lemma~28.23.5, Tag~0G41]{Stacks}; after removing torsion and taking the double dual it is a vector bundle $V$, because reflexive sheaves on a regular surface are locally free \cite[Lemma~31.13.15, Tag~0B3N]{Stacks}. The line bundle $L^\circ$ extends after taking the closure of a divisor. Finally, after twisting by a sufficiently large boundary divisor, $Q^\circ$ extends to
\[
 Q\in H^0\!\left(\overline T,\operatorname{Sym}^2V^\vee\otimes L\right).
\]
Multiplication of $Q$ by a boundary equation changes neither the conic nor the projective reflection on $T^\circ$. Let $\mathcal C$ be the irreducible component of $\{Q=0\}\subset\PP(V)$ dominating $\overline T$. It is birational to $\mathscr N$ and agrees with $\mathscr N^\circ$ over $T^\circ$.

For a section $v$ of a twist $V\otimes A$, the formula
\begin{equation}\label{eq:downstairs-reflection}
 R_v(w)=Q(v)w-2B_Q(w,v)v
\end{equation}
defines a rational involution whenever $Q(v)$ is not identically zero. Conjugation by the fixed birational comparison with $\mathscr N$ defines a right-fiber-preserving rational involution $\theta_v$ of $X_R$. Over $B^\circ$ it is
\[
 \theta_v=\pi\rho_v\pi^{-1},
\]
where $\rho_v$ is the reflection on $\mathscr N^\circ$.

\paragraph{Controlling the fibers outside $B^\circ$.}
The conic model has a direct geometric description over $B^\circ$. Along the boundary we prescribe sufficiently many terms of the reflection vector to ensure that all allowed reflections are regular on a fixed dense open subset of each remaining right fiber. A \emph{finite jet condition} means precisely such a congruence modulo a power of the ideal of the boundary.
We denote the extended projection by $\gamma:\overline T\to\PP_R^1$ and the set of omitted points by $B_{\partial}=\PP_R^1\setminus B^\circ$. The corresponding divisor on $\overline T$ is
\[
 D_{\partial}=\gamma^*\Bigl(\sum_{c\in B_{\partial}}[c]\Bigr).
\]
The boundary of $T^\circ$ lies over $B_{\partial}$: its open immersion into $\overline T|_{B^\circ}$ is proper, since $T^\circ$ is proper over $B^\circ$, and hence is also closed. Each surface $S_c=\varphi_R^{-1}(c)$ is integral.

\begin{lemma}\label{lem:boundaryjets}
There are fixed closed subsets $K_c\subset S_c$, each of dimension at most one, and an affine family of sections $v$ such that every birational member satisfies
\begin{equation}\label{eq:boundarycontainment}
 \Indet(\theta_v)\cap S_c\subset K_c
 \qquad(c\in B_{\partial}).
\end{equation}
The affine family is defined by a finite jet condition along $D_{\partial}$.
\end{lemma}
\begin{proof}
We start with a twist $A_0$ and a section $v_*$ whose reflection is birational. For $c\in B_{\partial}$, let $\eta_c$ be the generic point of $S_c$. The rational comparison $X_R\dashrightarrow\mathcal C$ is regular at $\eta_c$: the local ring $\OO_{X_R,\eta_c}$ is a discrete valuation ring, and properness of $\mathcal C$ extends the rational map over that DVR. Its parameter map to $\overline T$ is therefore regular there as well.

We embed $X_R$ in projective space and choose local frames for $V\otimes A_0$ near the comparison image of $\eta_c$. After clearing fixed denominators, we can write the coordinates of $\theta_v$ as
\[
 [H_0(v):\cdots:H_M(v)],
\]
where each $H_i(v)$ is polynomial in the components of $v$ and rational on $X_R$. Let $r_c$ be a local equation for $S_c$. Multiplying all coordinates by one fixed power of $r_c$ makes their coefficients regular at $\eta_c$. For the reference section $v_*$, the smallest order of its coordinates is
\[
 m_c=\min_i\operatorname{ord}_{S_c}H_i(v_*).
\]
After dividing the whole tuple by $r_c^{m_c}$, at least one coordinate has nonzero reduction on $S_c$.

If $n>m_c$ and $v$ agrees with $v_*$ modulo $r_c^n$, then every $H_i(v)-H_i(v_*)$ is divisible by $r_c^n$. After division by $r_c^{m_c}$ the two tuples therefore have the same reduction along $S_c$. Shrinking to a fixed neighborhood $U_c$ of $\eta_c$, the divided tuple is regular for every such $v$, and its common zero set on $U_c\cap S_c$ is contained in the fixed common zero set of the reduction of $r_c^{-m_c}H(v_*)$. Let $K_c$ be the union of $S_c\setminus U_c$ with the closure of that common zero set. It is a proper closed subset of the integral surface $S_c$, hence has dimension at most one. Outside $K_c$ the tuple defines a morphism to projective space, whose image lies in $X_R$ because this holds on a dense open subset. This proves \eqref{eq:boundarycontainment} for the jet condition at $c$.

There are only finitely many points $c\in B_{\partial}$. Taking one integer $n$ larger than all $m_c$ imposes these conditions simultaneously along $nD_{\partial}$.
Indeed, on a fixed neighborhood of $\eta_c$ the parameter comparison $\kappa$ is a morphism and $\gamma\circ\kappa=\varphi_R$. The pullback of the ideal of $nD_{\partial}$ is therefore contained in $(r_c^n)$, so the prescribed jets imply the congruence used in the coordinate calculation.
\end{proof}

\paragraph{Retaining freedom to choose the reflection.}
The boundary conditions must leave enough freedom to prescribe values at two distinct points of $T^\circ$. This independence will be used to control the left surface through the moving point $y$. With $n$ fixed, we increase the twist to $A=A_0\otimes\mathcal H^k$, where $\mathcal H$ is ample on $\overline T$. For sufficiently large $k$ there is a section $a_k\in H^0(\overline T,\mathcal H^k)$ nonzero at the parameter $t_0\in T_b$ of the conic through $x$ and at each generic comparison image above $S_c$. With this choice the sections satisfying the prescribed boundary condition form the affine space
\begin{equation}\label{eq:jetfamily}
 P_0=a_kv_*+H^0\bigl(\overline T,V\otimes A(-nD_{\partial})\bigr).
\end{equation}
Near those generic images, $a_k$ is a unit and $v/a_k$ agrees with $v_*$ to order $n$. Scaling a vector leaves its projective reflection unchanged. Thus \eqref{eq:boundarycontainment} holds for this family, after including the zeros of $a_k$ in the fixed complements $S_c\setminus U_c$.

For sufficiently large $k$, the variable part
\[
 H^0(\overline T,V\otimes A(-nD_{\partial}))
\]
allows arbitrary values at any two distinct points of $T^\circ$. For a bound on $k$ uniform over all pairs, we apply relative Serre vanishing to the ideal sheaf of the universal length-two subscheme over $\operatorname{Hilb}^2(\overline T)$, tensored with the pullback of the bundle in the preceding display. Evaluation is then surjective for every pair once $k$ is large. The restriction map to the complete central curve $T_b$ is also surjective for large $k$, by Serre vanishing applied to the exact sequence obtained by twisting with $-T_b$. Since $D_{\partial}$ is disjoint from $T_b$, the prescribed boundary jets impose no condition on that restriction.

\paragraph{Conditions on the central fiber.}\label{app:centralconditions}
The surjectivity of restriction to $T_b$ allows the choices from Section~\ref{sec:relocation}. There is a nonempty open subset $P\subset P_0$ on which the following properties hold:
\begin{enumerate}
\item $Q(v)$ is not identically zero and $Q(v(t_0))\ne0$;
\item $\theta_v$ is an isomorphism between neighborhoods of $x$ and $y(v)=\theta_v(x)$;
\item $y(v)\notin C_R$, and its image in $X\setminus\Sigma$ lies on no line in either of the two corresponding del Pezzo surfaces in $\PP^4$;
\item $\rk d(\varphi_R,\psi_L)_{y(v)}=2$;
\item $\Indet(\theta_v^{-1})\cap S_b$ is finite and misses all sixteen marked points.
\end{enumerate}
The parameter $t_0$ and the conic $(h^\circ)^{-1}(t_0)$ through the lift of $x$ are fixed. As the anisotropic value $v(t_0)$ varies, the reflected point ranges over an open subset of this conic. We choose that point away from the auxiliary exceptional divisor, where $\pi$ is an isomorphism. We also exclude the lift of $x$, the only lift of a marked point on this conic by Lemma~\ref{lem:auxiliarysection}(1). The line and critical-locus conditions of Subsection~\ref{app:oppositekernel} exclude only finitely many further points. Thus (1)--(4) hold for a nonempty open set of values of $v(t_0)$.

For (5), we require $v|_{T_b}$ to be nowhere zero and anisotropic at the parameters of the singular conics and of the conics through the lifts of the $s_j(b)$. At a zero of $Q(v)$, the conic is smooth and the reflection formula has a single possible base point. The comparison $\pi^{-1}$ adds only $q(b)$ and is regular near each $s_j(b)$, and the reflection is regular on the complete conic through its lift. Hence $\theta_v$ is regular at all sixteen marked points and has only finitely many possible indeterminacy points in $S_b$. The interpolation in Section~\ref{sec:reflectioninterpolation} imposes these conditions simultaneously. Nonvanishing on the complete curve and the stated finite-parameter exclusions are open conditions on sections, so they define a nonempty open subset $P$ with all five properties.

After shrinking $P$, the map
\[
 \lambda:P\to\PP_L^1,\qquad \lambda(v)=\psi_L(y(v)),
\]
is a morphism. It depends only on $v(t_0)$. As $v(t_0)$ varies, the image of $x$ moves through an open subset of the normalization of $D_0$. The left parameter is nonconstant on $D_0$, see Lemma~\ref{lem:opposite-nonconstant}. Hence $\lambda$ is nonconstant as a function of $v(t_0)$.

\Needspace{9\baselineskip}
\subsection{The left surface through the image point}
It remains to show that the reflection can be chosen so that its indeterminacy does not contain a curve in $S_L(y(v))$. The point $y(v)$, and hence this left surface, vary with $v(t_0)$. At a point of another conic we can vary the reflection independently, using the freedom of two-point evaluation established above.

\begin{lemma}\label{lem:oppositefiber}
For general $v\in P$, the intersection
\[
 \Indet(\theta_v^{-1})\cap S_L(y(v))
\]
is finite.
\end{lemma}
\begin{proof}
We divide the proof according to the location of a possible indeterminacy point. On the open conic model we count dimensions while varying $v$. The remaining points belong either to the fixed curves obtained from the boundary construction or to $C_R$, the union of the exceptional curves $C_p^R$.

\emph{Indeterminacy on the open conic model.}

For $z\in X_R|_{B^\circ}\setminus q(B^\circ)$, let $[w(z)]$ be its lift to $\mathscr N^\circ$ and $t(z)=h^\circ(w(z))$. The possible zeros of the reflection formula belong to the incidence set
\[
 \mathcal E=\{(v,z):R_{v(t(z))}(w(z))=0\}.
\]
Outside $\mathcal E$ the reflection formula defines a morphism into $\mathscr N^\circ$, followed by the morphism $\pi$. Hence
\begin{equation}\label{eq:reflectionbadset}
 \Indet(\theta_v)\cap
 (X_R|_{B^\circ}\setminus q(B^\circ))\subset\mathcal E_v.
\end{equation}

Let $d=\dim P$ and let $\Delta\subset T^\circ$ be the conic discriminant. On a smooth conic, for a fixed point $[w]$, the equation $R_v(w)=0$ holds if and only if $v\in\CC w$. Indeed, the reflection associated with an anisotropic vector is invertible. For nonzero isotropic $v$, the equation requires $B_Q(w,v)=0$, and the tangent line to a smooth conic meets it only at $[v]$. This imposes codimension two on the three-dimensional space of values $v(t)$.

On a singular conic the equation imposes codimension at least one: it is a proper closed condition, since every nonzero quadratic form admits an anisotropic vector, whose reflection is invertible. Surjectivity of point evaluation identifies the value $v(t)$ with a free three-dimensional parameter. For points over smooth conics, $z$ varies in a three-dimensional set and the condition on $v$ has codimension two. The corresponding incidence therefore has dimension at most $d+3-2=d+1$. Over singular conics, the set of possible $z$ has dimension at most two, and the condition on $v$ has codimension at least one. The dimension bound in this case is again $d+2-1=d+1$. Consequently $\dim\mathcal E\le d+1$.

\emph{Requiring the point to lie on the moving left surface.}
For $z\notin C_R$ with $t(z)\ne t_0$, impose also
\[
 \lambda(v)=\psi_L(z).
\]
The two-point evaluation property makes $v(t(z))$ and $v(t_0)$ independent coordinates on the parameter space. The reflection equation constrains only the first value, while $\lambda(v)$ depends only on the second. Since $\lambda$ is nonconstant as $v(t_0)$ varies, the equation $\lambda(v)=\psi_L(z)$ cuts the permissible values at $t_0$ by at least one further condition. The dimension bound thus decreases from $d+1$ to $d$ in both cases. In other words, the incidence of pairs $(v,z)$ for which the reflection formula vanishes and $z$ lies on the left surface through $y(v)$ has dimension at most $\dim P$.
There are no points of $\mathcal E$ over $t_0$, because $Q(v(t_0))\ne0$. By the fiber-dimension theorem, for general $v\in P$,
\begin{equation}\label{eq:finiteopenincidence}
 \{z\in\mathcal E_v\setminus C_R:\psi_L(z)=\lambda(v)\}
 \quad\text{is finite}.
\end{equation}

\emph{The auxiliary section and the omitted right fibers.}
The remaining possible indeterminacy points outside $C_R$ lie in the fixed set
\[
 K=q(\PP_R^1)\cup\bigcup_{c\in B_{\partial}}K_c,
\]
which has dimension at most one. On an irreducible curve of $K$ not contained in $C_R$, the rational map $\psi_L$ is either nonconstant, with finite level sets on its domain, or constant. Since $\lambda$ is nonconstant, the parameters $v$ for which $\lambda(v)$ avoids the finitely many constant values form a nonempty open subset of $P$.

\emph{The exceptional curves $C_p^R$.}
Finally, $\Indet(\theta_v)$ has codimension at least two. Indeed, $X_R$ is smooth, hence normal, and a rational map from a normal variety to a projective variety extends across every codimension-one point by the valuative criterion for properness. Since $X_R$ is a threefold, $\dim\Indet(\theta_v)\le1$. It contains no curve $C_p^R$: by item~(5), $\theta_v$ is regular at a point of each such section. Its intersection with $C_R$ is therefore finite. Away from $C_R$, the closed surface $S_L(y(v))$ agrees with the level set of $\psi_L$. Equations~\eqref{eq:boundarycontainment}, \eqref{eq:reflectionbadset}, and~\eqref{eq:finiteopenincidence} show that its intersection with $\Indet(\theta_v)$ is finite on each of the regions just considered. Since $\theta_v$ is an involution, this is also the assertion for $\theta_v^{-1}$.
\end{proof}

\begin{proof}[Proof of Theorem~\ref{thm:smallrelocation}]
For a general $v\in P$, we take $\theta=\theta_v$. By conditions (1)--(4) in Section~\ref{app:centralconditions}, $\theta$ is a local isomorphism at $x$, its image $y$ lies in the required locus, and the two differentials at $y$ are independent. By (5), the indeterminacy on the right fiber is finite and avoids its marked points. The intersection with the left surface is finite, see Lemma~\ref{lem:oppositefiber}. Finally, $\theta$ is regular at a point of each $C_p^R$, and its indeterminacy has dimension at most one, so its intersection with $C_R$ is finite. Since $\theta$ is an involution, these statements hold for its inverse as well.
\end{proof}

\subsection{Constructing the three independent spray directions}
We can now construct the sprays. At the image point $y$, we use two independent directions along the right fibration and one direction along the left fibration outside their span. The finiteness assertions just proved allow the complete curves for all three sprays to be chosen disjoint from inverse indeterminacy.

\begin{proposition}\label{prop:downstairs-sprays}
Every point of $X_R$ admits a dominating finite family of local algebraic sprays taking values in $X_R$.
\end{proposition}
\begin{proof}
If $x\notin C_R$, let $\widetilde x\in Z$ be its unique lift. A locally dominating family on $Z$ at $\widetilde x$ may be based on a neighborhood disjoint from the exceptional divisors. The evaluation maps composed with $\beta_R$ form a dominating local family at $x$, because $d\beta_R$ is an isomorphism at $\widetilde x$.

For $x\in C_R$, Theorem~\ref{thm:smallrelocation} supplies $\theta$ and $y=\theta(x)$. We denote the right parameter of $y$ by $b=\varphi_R(y)$ and its unique lift to $Z$ by $\widetilde y$. For the inverse indeterminacy locus $J=\Indet(\theta^{-1})$, we write
\[
 \widetilde J=\beta_R^{-1}(J).
\]
The intersection of $\widetilde J$ with the right $g_R$-fiber through $\widetilde y$ is finite. We have $J\cap S_b\cap C_R=\varnothing$, see Theorem~\ref{thm:smallrelocation}, so its inverse image contains no exceptional curve in this right fiber. The intersection with the left $g_L$-fiber is also finite: away from $C_R$ this follows from Theorem~\ref{thm:smallrelocation}; for $x'\in J\cap C_p^R$, the exceptional fiber $\beta_R^{-1}(x')\subset E_p$ meets a fixed left $g_L$-fiber in one point.

The finite forbidden sets in these two fibers allow us to apply Lemmas~\ref{lem:ordinarycurves} and~\ref{lem:relativeordinary}. The two branches of one right-fiber nodal curve determine right-relative sprays at $\widetilde y$ with independent velocities
\[
 v_1,v_2\in\ker(dg_R)_{\widetilde y}.
\]
For a nodal curve in the left fiber, at least one of its two independent branch directions lies outside the right vertical plane, since the two planes are distinct. The corresponding left-relative spray has velocity
\[
 v_3\in\ker(dg_L)_{\widetilde y}\setminus\ker(dg_R)_{\widetilde y}.
\]
The nodal curves are chosen on the del Pezzo surfaces to avoid the projections of $\widetilde J$ and the marked points. Their lifts therefore avoid $\widetilde J$. Lemma~\ref{lem:relativeordinary} extends them to proper families whose complete images still avoid $\widetilde J$ after shrinking the base neighborhoods. Restricting these families to the complements of their infinity sections defines the three sprays above. The hypotheses of the lemmas hold because $y$ lies on no line in either corresponding del Pezzo surface, and the forbidden sets are finite. Since $\rk d(g_R,g_L)_{\widetilde y}=2$, the three velocities span $T_{\widetilde y}Z$.

We shrink the base neighborhood of $\widetilde y$ to avoid the exceptional divisors. Via $\beta_R$ we identify it with a neighborhood $U_y\subset X_R$ and compose the three evaluation maps with $\beta_R$. The resulting sprays are
\[
 s_i:U_y\times\AA^1\longrightarrow X_R\setminus J
\]
whose velocities span $T_yX_R$. Lemma~\ref{lem:conjugation} transports this dominating family to $x$.
\end{proof}

\begin{proof}[Proof of Theorem~\ref{thm:smallres}]
The variety $X_R$ is algebraically elliptic, see Proposition~\ref{prop:downstairs-sprays} and Lemma~\ref{lem:localglobal}. The hypotheses used to construct the sprays are symmetric in the two pencils, so the same proof applies to $X_L$ with right and left interchanged. Condition~(iii) of Notation~\ref{nota-U} enters only the proof of nonrationality for $X_R$ in Corollary~\ref{cor:nonrationality}. Since $X_L$ is birational to $X_R$, it too is nonrational. This completes the proof of Theorem~\ref{thm:main}.
\end{proof}

\appendix
\section{Surface geometry for \texorpdfstring{$C_*$}{C*}}\label{app:geometry}
We check the conditions defining $\mathcal U$ for the matrix $C_*$ in \eqref{eq:witnessC}. The calculations concern three geometric properties: the singularities of the del Pezzo surfaces, the intersections with coordinate subspaces, and the existence of a fiber containing no line through any marked point.

\subsection{The singular fibers of the two fibrations}
At a homogeneous parameter $[u:v]$, the right coefficients are
\[
 a_i=uc_{0i}+vc_{1i},\qquad b_i=uc_{2i}+vc_{3i},\qquad
 D^R_{ij}=a_ib_j-a_jb_i.
\]
The left coefficients use the row pairs $(0,2)$ and $(1,3)$. The determinant of the matrix $\mathsf C_i$ tests whether the $i$th coefficient column can vanish for some $[u:v]$. For $i=0,\ldots,4$, these determinants are
\[
 2,\quad -2,\quad -2,\quad -3,\quad 1.
\]
Thus no coefficient column vanishes. The following table lists the minors that determine the singular parameters. For each pair $ij$, the right and left binary quadratic forms have the same discriminant, shown in the last column.
\begin{center}
\small
\begin{tabular}{c@{\quad}llr}
\toprule
$ij$&$D^R_{ij}(u,v)$&$D^L_{ij}(u,v)$&Discriminant\\\midrule
$01$&$6u^2+14uv+6v^2$&$3u^2-4uv-3v^2$&$52$\\
$02$&$-4u^2-6uv+v^2$&$-2u^2+6uv+2v^2$&$52$\\
$03$&$2u^2-3uv-5v^2$&$-8u^2-9uv-v^2$&$49$\\
$04$&$-4u^2-5uv-2v^2$&$7u^2+7uv+2v^2$&$-7$\\
$12$&$2u^2+3uv+2v^2$&$u^2-3uv+4v^2$&$-7$\\
$13$&$-7u^2+uv+2v^2$&$u^2-5uv-8v^2$&$57$\\
$14$&$11u^2+10uv+2v^2$&$u^2+4uv+v^2$&$12$\\
$23$&$4u^2-3uv+2v^2$&$2u^2+3uv+4v^2$&$-23$\\
$24$&$-6u^2+v^2$&$-3u^2+2v^2$&$24$\\
$34$&$-u^2-5uv-v^2$&$-5u^2-11uv-5v^2$&$21$\\\bottomrule
\end{tabular}
\end{center}
All ten forms in each column are squarefree, and no two have a common projective zero. Thus either pencil has twenty distinct singular parameters, with exactly one proportional column pair at each.

\subsection{The coefficient and coordinate conditions}
The first four-column determinant of $C_*$ is $-45$, and its rows sum to zero. The ten $2\times2$ determinants comparing the coefficient forms $uc_{0i}+vc_{1i}$ are
\[
 3,-2,-8,7,1,1,1,2,-3,-5
\]
in the order $01,02,03,04,12,13,14,23,24,34$; hence the five roots used in condition~(iii) are distinct. In the squared coordinates $z_i=x_i^2$, the quartic form is $f_C(x)=z^TNz$. We use this expression to check its intersections with the coordinate subspaces. If $\mathbf c_k=(c_{k0},\ldots,c_{k4})$ denotes the $k$th row of $C_*$, then
\[
 N=\frac12\bigl(\mathbf c_0^T\mathbf c_3+\mathbf c_3^T\mathbf c_0-\mathbf c_1^T\mathbf c_2-\mathbf c_2^T\mathbf c_1\bigr).
\]
Every principal minor of $N$ of size one, two, or three is nonzero. This proves the coordinate-stratum smoothness assertions used in Section~\ref{sec:genericity}. The ten coordinate-line quartics in Appendix~\ref{app:rankzero} are squarefree with nonzero leading and constant coefficients.

\subsection{A fiber with no line through a marked point}
The finiteness argument in Lemma~\ref{lem:linefiniteness} requires one fiber of each fibration that contains no line through any $p\in\Sigma$. At the parameter $[1:1]$ all ten surface minors are nonzero in either pencil. The line polynomial \eqref{eq:linepolynomial} evaluated at a node $p\in\Sigma$ has value $1530720$ in the right pencil and $1027728$ in the left pencil, so this fiber contains no line through a node $p\in\Sigma$.

The exact calculations used here are described in Appendix~\ref{app:reproduction}.

\section{Node involutions on the rank-one locus for \texorpdfstring{$C_*$}{C*}}\label{app:rankone}
For $C_*$, we must show that every point with exactly two zero coordinates admits a suitable node in the sense of Definition~\ref{def:suitable-node}. Such points occur in ten coordinate planes. On each plane, we shall express the required geometric conditions by polynomials and show that the required polynomials are nonzero for at least one choice of signs. The sign changes then turn this choice into a choice of node.

\subsection{An equation for the union of the sixteen lines}
We first need a polynomial that detects whether a point of a smooth degree-four del Pezzo surface lies on one of its lines. We work on a smooth diagonal intersection
\[
 S=\left\{\sum_{i=0}^4a_iw_i^2=\sum_{i=0}^4b_iw_i^2=0\right\},
 \qquad D_{ij}=a_ib_j-a_jb_i.
\]
All ten $D_{ij}$ are nonzero. The coefficients we need are
\[
 \alpha=-D_{04}D_{14}D_{23},\qquad
 \beta=D_{03}D_{13}D_{24},\qquad
 \gamma=-D_{02}D_{12}D_{34}.
\]
In the squared coordinates $z_i=w_i^2$, the line test uses the polynomial
\begin{equation}\label{eq:linepolynomial}
 \Lambda=\gamma^2z_2^2+\beta^2z_3^2+\alpha^2z_4^2
 -2\beta\gamma z_2z_3-2\alpha\gamma z_2z_4-2\alpha\beta z_3z_4.
\end{equation}
\begin{lemma}\label{lem:line-polynomial-test}
On $S$, the zero locus of $\Lambda$ is the union of its sixteen lines.
\end{lemma}
\begin{proof}
We choose a basis of the pencil in which all $a_i$ are nonzero. The ratios $\lambda_i=b_i/a_i$ are then pairwise distinct. If $f(T)=\prod_i(T-\lambda_i)$, then
\[
 w_i^2=\frac{(s+\lambda_i t)^2}{a_if'(\lambda_i)}
\]
are the squared images of sixteen lines obtained by changing coordinate signs. The two quadric identities follow from
\[
 \sum_i\frac{\lambda_i^j}{f'(\lambda_i)}=0,\qquad 0\le j\le3.
\]
After substitution of these squared coordinates into \eqref{eq:linepolynomial}, the expression vanishes after cancellation of the factors $\lambda_i-\lambda_j$. Thus $\Lambda$ vanishes on all sixteen lines. To see that there are no other components, the two linear equations in the $z_i$ determine $z_0,z_1$ from $z_2,z_3,z_4$, since $D_{01}\ne0$. The nonzero polynomial in \eqref{eq:linepolynomial} therefore does not vanish identically on the squared-coordinate plane, nor on $S$. Its zero divisor on $S$ has class $4H$ and degree $16$, where $H=\OO_S(1)$. The sixteen distinct lines already have total degree $16$, so they exhaust its support. See also \cite[Section~2]{Skoro} for the diagonal description of these lines.
\end{proof}

\subsection{Parametrizing a critical component by squared coordinates}
For $I=\{a,b,c\}$ with $a<b<c$, we work on the coordinate plane $x_i=0$ for $i\notin I$. We write $Z_i=x_i^2$. The two equations of the right del Pezzo surface are linear in $Z_a,Z_b,Z_c$. Their projective solution is
\begin{equation}\label{eq:quotientappendix}
 [Z_a:Z_b:Z_c]=[D^R_{bc}:-D^R_{ac}:D^R_{ab}].
\end{equation}
For $C_*$, the three quadratics are linearly independent for every $I$; their coefficient determinants are respectively
\[
 50,314,-46,234,-106,-79,116,71,-81,101
\]
for $I=012,013,014,023,024,034,123,124,134,234$. Thus the map in \eqref{eq:quotientappendix} identifies $\PP^1$ with the conic of squared coordinates of this critical component. Passing back from the $Z_i$ to the $x_i$ requires choices of square roots. We keep all of these choices in one algebra. In the chart $[u:v]=[t:1]$, work in
\begin{equation}\label{eq:rank8}
 R_I=\QQ[t,\xi_a,\xi_b,\xi_c]/
 (\xi_a^2-Z_a,\xi_b^2-Z_b,\xi_c^2-Z_c).
\end{equation}
Every element of $R_I$ has a unique expression over $\QQ[t]$ in the eight monomials $\xi_a^{e_a}\xi_b^{e_b}\xi_c^{e_c}$, where each $e_i$ is $0$ or $1$. At a parameter with $Z_aZ_bZ_c\ne0$, these expressions can be evaluated at the eight choices of square roots. Simultaneous reversal of all three signs represents the same projective point, leaving four projective choices.

\paragraph{The involution and its image.}
We take $x_i=\xi_i$ on $I$ and $x_i=0$ off $I$. At $p=(1,1,1,1,1)$, let $L$ be the matrix formed from the four entries of $2Cx$ in the same order as $M$. With $m=M(x)$ we have
\[
 A=\det L,\qquad
 B=L_{00}m_{11}+m_{00}L_{11}-L_{01}m_{10}-m_{01}L_{10},\qquad y=Ax-Bp.
\]
Polynomial kernel coordinates for $M(y)$ are obtained from
\[
 r_y=\operatorname{adj}(L)c(t),\qquad
 \ell_y=(1,-t)\operatorname{adj}(L),
\]
where $c(t)$ is the first column of $M(x)$, with entries in $R_I$. The identity $M(x)(t,1)^T=0$ on the quotient conic implies $M(x)=c(t)(1,-t)$. The verification checks these kernel identities and $f_C(y)=0$ in \eqref{eq:rank8}.

\paragraph{Expressing the geometric conditions by one polynomial.}
Let $\Delta_\epsilon(y)$ be the product of the ten coefficient minors of the $\epsilon$-surface at the parameter computed from the kernel formulas. We write $\Lambda_\epsilon(p)$ and $\Lambda_\epsilon(y)$ for the values of \eqref{eq:linepolynomial} on that surface and form the product
\begin{equation}\label{eq:goodpoly}
 P_I=A B^2\left(\prod_{i=0}^4y_i^2\right)
 \prod_{\epsilon=R,L}\bigl(\Delta_\epsilon(y)
          \Lambda_\epsilon(p)\Lambda_\epsilon(y)\bigr).
\end{equation}
If $P_I\ne0$, then $A\ne0$ and every coordinate of $y$ is nonzero. The two computed parameters define smooth fibers of the two del Pezzo fibrations containing no line through $p$ or $y$. We also verify below that $M(y)\ne0$, so these are the parameters of the right and left fibers through the nonnodal point $y$.

\paragraph{Comparing the four sign choices.}
In $R_I$ we reduce powers of $\xi_i$ using $\xi_i^2=Z_i$. The polynomial $P_I$ is invariant under simultaneous reversal of all three signs. Consequently only the four monomials of even total degree remain:
\begin{equation}\label{eq:characters}
 P_I=P_0+\xi_a\xi_bP_{ab}+\xi_a\xi_cP_{ac}+\xi_b\xi_cP_{bc},
\end{equation}
with degrees $238,236,236,236$. When $Z_aZ_bZ_c\ne0$, there are four values of $P_I$ under the projective sign choices. These values are linear combinations of $P_0$, $\xi_a\xi_bP_{ab}$, $\xi_a\xi_cP_{ac}$, and $\xi_b\xi_cP_{bc}$ with an invertible matrix of signs. Hence all four values vanish if and only if all four coefficient polynomials vanish. This replaces the choice among four projective points by a common-zero test for four polynomials in one variable. If $\delta$ is a coordinate sign change, then
\[
 \tau_{\delta p,C}(x)=\delta\,\tau_{p,C}(\delta x).
\]
Thus changing the coordinate signs of the source for the fixed center $p=(1,1,1,1,1)$ is equivalent, for the geometric conditions defining \(\mathscr G_p\), to varying the center $p\in\Sigma$ while keeping the source point fixed.

\subsection{Finding a sign choice for which the node is suitable}
We first check that the four coefficient polynomials have no common zero away from $Z_aZ_bZ_c=0$. For polynomials in one variable this is equivalent to checking the roots of their greatest common divisor. For $C_*$, the computed gcd is
\begin{equation}\label{eq:unifiedgcd}
 \gcd(P_0,P_{ab},P_{ac},P_{bc})=Z_aZ_bZ_c
\end{equation}
up to a nonzero rational scalar. This excludes simultaneous vanishing away from the coordinate boundary in the chart $[u:v]=[t:1]$. At the remaining parameter $[u:v]=[1:0]$, we use the second column of $M(x)$ and the row factor $(0,1)$ in the formulas for the kernels of $M(y)$. The reciprocal-chart calculation verifies $P_I\ne0$ and $M(y)\ne0$ simultaneously for all four projective points above these squared coordinates on each of the ten supports. It uses nonzero norms in the constant rank-eight algebra: the norm of $P_I$ is nonzero, and at least one entry of $M(y)$ has nonzero norm. Thus the omitted parameter is covered as well.

\paragraph{Excluding an image that is itself a node.}
The condition $P_I\ne0$ does not alone exclude $M(y)=0$. We must ensure that the image for a sign choice satisfying the previous conditions is also nonnodal. On the affine chart, the calculation therefore also uses the four products $P_I M_{rs}(y)$. We expand each product in the four monomials of \eqref{eq:characters}. The greatest common divisor of all their coefficient polynomials again has zeros only where $Z_aZ_bZ_c=0$. The values at the four sign choices are obtained from these coefficients by the same invertible matrix of signs. Therefore, at a point off the coordinate boundary, it cannot happen that every sign for which $P_I$ is nonzero has $M(y)=0$: if it did, all coefficient polynomials of all products $P_IM_{rs}(y)$ would vanish there. Hence there is a sign choice for which both $P_I\ne0$ and $M(y)\ne0$. By sign symmetry there is consequently a center with $(C_*,x)\in\mathscr G_p$.

The calculation also verifies the identities for the kernels of $M(y)$ and $A_p(y)=A_p(x)^3$. The computational verification is described in Appendix~\ref{app:reproduction}.

\section{Node involutions on the rank-zero locus for \texorpdfstring{$C_*$}{C*}}\label{app:rankzero}
We now consider the points with exactly three zero coordinates. There are forty such points for $C_*$, four on each coordinate line. We shall check the stronger assertion that every one of the sixteen nodes is suitable for each of them. First we check the image of the involution by polynomial calculations. We then use the description of the lines on a two-nodal del Pezzo surface to check the remaining condition geometrically.

\subsection{Locating the forty critical points}
On the coordinate line indexed by $i<j$, we write $x=e_i+t e_j$. Up to nonzero rational scalars, the ten restrictions $f_{ij}(t)=f_C(e_i+t e_j)$ are
\begin{center}
\begin{tabular}{c@{\quad}lc@{\qquad}c@{\quad}l}
\toprule
Support&$f_{ij}(t)$& &Support&$f_{ij}(t)$\\\midrule
$01$&$t^4+3t^2-1$&&$12$&$2t^4-3t^2+2$\\
$02$&$t^4-3t^2-1$&&$13$&$3t^4-9t^2+2$\\
$03$&$(t^2+2)(3t^2-1)$&&$14$&$t^4-2t^2-2$\\
$04$&$t^4+t^2+2$&&$23$&$3t^4+t^2+2$\\
$24$&$t^4-4t^2-2$&&$34$&$t^4+3t^2-3$\\\bottomrule
\end{tabular}
\end{center}
All are squarefree with nonzero constant and leading terms,
so they account for forty reduced rank-zero points.

\subsection{Checking that the image is a nonnodal point of \texorpdfstring{$\Omega$}{Omega}}
For each pair $i<j$ and each $p\in\Sigma$, let $x=e_i+t e_j$ and let
\[
 y=A_p(x)x-B_p(x)p.
\]
The computed gcd is
\begin{equation}\label{eq:rankzerogcd}
 \gcd\!\left(f_{ij},A_p(x)B_p(x)\prod_{k=0}^4y_k\right)=1
\end{equation}
for all ten pairs $i<j$ and all sixteen nodes $p$, for a total of $160$ polynomial tests.
Each test treats all four roots of the corresponding $f_{ij}$ at once.
The computation also verifies that no root of $f_{ij}$ makes all four entries of $M(y)$ vanish.
Thus, for every root of $f_{ij}$ and every node $p$, we have $A_p(x)\ne0$, and $y$ is a nonnodal point in the coordinate torus. Its lift belongs to $\Omega$ by Lemma~\ref{lem:critical-coordinates}. By \eqref{eq:Aidentity},
the node involution is an isomorphism between neighborhoods of $x$ and $y=\tau_p(x)$.

\subsection{Excluding lines through the center and the image}
It remains to show that neither del Pezzo surface through $y$
contains a line through $p$ or $y$. The $k$th column of $C$ corresponds to the matrix
\[
 \mathsf C_k=\begin{pmatrix}c_{0k}&c_{1k}\\c_{2k}&c_{3k}\end{pmatrix}.
\]
Since the rows of $C$ sum to zero, $\sum_k\mathsf C_k=0$. The ten values
\begin{center}
\begin{tabular}{c@{\quad}rrrrrrrrrr}
$ij$&01&02&03&04&12&13&14&23&24&34\\\midrule
$f_C(e_i+e_j)$&$-6$&$6$&$-6$&$4$&$-1$&$4$&$-3$&$-6$&$-5$&$1$
\end{tabular}
\end{center}
are all nonzero; equivalently, every $\mathsf C_i+\mathsf C_j$ is invertible.

Let $r$ be a right-kernel vector of $M(y)$. The parameter of the right surface through $y$ is $[r]$. The $i$th and $j$th coefficient columns of that surface are the vectors $\mathsf C_i r$ and $\mathsf C_j r$. Since $y_k^2=B_p(x)^2$ for $k\notin\{i,j\}$ and $\sum_k\mathsf C_k=0$, the equation $M(y)r=0$ reads
\[
 (y_i^2-B_p(x)^2)\mathsf C_i r+(y_j^2-B_p(x)^2)\mathsf C_j r=0.
\]
The two scalar coefficients cannot both vanish, since that would imply $M(y)=0$. Neither vector is zero, since every $\mathsf C_k$ is invertible by condition (i). Thus the two columns are proportional; by condition (ii) they are the unique proportional pair in that singular fiber. A line through the node $p$ would, by \eqref{eq:singularlinecriterion}, imply
\[
 (\mathsf C_i+\mathsf C_j)r=0,
\]
which is impossible. The same kernel equation expresses the line criterion at $y$ as
\[
 (y_i^2\mathsf C_i+y_j^2\mathsf C_j)r
   =B_p(x)^2(\mathsf C_i+\mathsf C_j)r\ne0.
\]
The same argument with a left-kernel row vector treats the left surface. Neither of the two del Pezzo surfaces through $y$ contains a line passing through $p$ or $y$. Thus Proposition~\ref{prop:ordinaryavoid} applies at every point of $R_0$ for every center $p\in\Sigma$.

\section{Exact computational verification}\label{app:reproduction}
The computations in the preceding three appendices use exact arithmetic over $\QQ$. The verifier performs the following checks.
\begin{enumerate}
\item For membership in $\mathcal U$, it computes the coefficient determinants and the binary quadratic minors of Appendix~\ref{app:geometry}. It checks their nonvanishing, discriminants, and absence of pairwise common projective zeros, as well as the principal minors of $N$ and the line tests at $[1:1]$.
\item On each of the ten components of $R_1$, it reduces the polynomials in the algebra \eqref{eq:rank8} and computes the greatest common divisors of their coefficient polynomials. In particular, the products $P_I M_{rs}(y)$ test whether a single sign choice gives both $P_I\ne0$ and $M(y)\ne0$. Together with sign symmetry, these tests show that a suitable node exists at each point of the affine chart on this component. The omitted parameter $[1:0]$ is checked separately, as explained in Appendix~\ref{app:rankone}.
\item On each coordinate line, it computes \eqref{eq:rankzerogcd} for each of the sixteen nodes, and checks that the entries of $M(y)$ do not vanish simultaneously at a root of $f_{ij}$. Together with the geometric line test in Appendix~\ref{app:rankzero}, these calculations treat all forty points of $R_0$.
\end{enumerate}
The verification also checks the formulas for the right and left kernels of $M(y)$ and the identities $f_C(y)=0$ and $A_p(y)=A_p(x)^3$ used in these calculations.

The exact-arithmetic verifier for the computational assertions in Appendices~\ref{app:geometry}--\ref{app:rankzero}, together with dependency specifications, execution instructions, and reference output, is available on Zenodo (version~1.0.0, DOI \href{https://doi.org/10.5281/zenodo.24345191}{\texttt{10.5281/zenodo.24345191}}). It verifies that the displayed coefficient matrix $C_*$ belongs to $\mathcal U_0$. The reference output includes all ten checks at $[u:v]=[1:0]$ and the number of nodes satisfying conditions (a)--(c) for each coordinate line.

\end{document}